\documentclass{gtart}

\def\ifplaintex{\expandafter\ifx\csname documentclass\endcsname\relax}

\expandafter\ifx\csname epsfbox\endcsname\relax\input epsf\fi

\ifplaintex 
\else
\fi

\def\gt{{\mathsurround=0pt\it $\cal G\mskip-2mu$eometry \&\ 
$\cal T\!\!$opology}}        

\def\gtp{{\mathsurround=0pt\it $\cal G\mskip-2mu$eometry \&\ 
$\cal T\!\!$opology $\cal P\!$ublications}}  

\def\lognumber#1{\def\thelognumber{#1}}
\def\volumenumber#1{\def\thevolumenumber{#1}}
\def\papernumber#1{\def\thepapernumber{#1}}
\def\volumeyear#1{\def\thevolumeyear{#1}}

\def\pagenumbers#1#2{\def\startpage{#1}\def\finishpage{#2}}
\def\published#1{\def\publishdate{#1}}
\def\proposed#1{\def\theproposer{#1}}
\def\seconded#1{\def\theseconders{#1}}
\def\received#1{\def\receiveddate{#1}}

\def\accepted#1{\def\accepteddate{#1}}
\def\asciititle#1{\def\theasciititle{#1}}
\def\covertitle#1{\def\thecovertitle{#1}}

\long\def\asciiabstract#1{\long\def\theasciiabstract{#1}}

\let\\\par\let\thevolumenumber\relax\let\thepapernumber\relax
\let\thevolumeyear\relax\let\thesamplenumber\relax\let\startpage\relax
\let\finishpage\relax\let\publishdate\relax\let\receiveddate\relax
\let\reviseddate\relax\let\accepteddate\relax\let\theasciititle\relax
\let\thecovertitle\relax\let\theasciiauthors\relax
\let\theasciiabstract\relax
\let\theasciiemail\relax\let\theshortauthors\relax\let\theshorttitle\relax

\long\def\maketitlep{   

\count0=\startpage

\gt\hfill      
\hbox to 77pt{\vbox to 0pt{\vglue -15pt\epsfbox{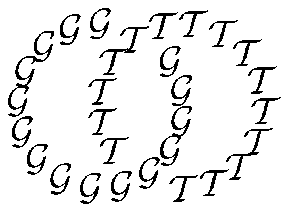}\vss}\hss}
\break
{\small\ifx\thesamplenumber\relax 
Volume \else Sample
\fi\thevolumenumber\ (\thevolumeyear)
\startpage--\finishpage\nl
Published: \publishdate}
\vglue 0.5truein plus 0.4fil minus 0.1truein

{\parskip=0pt\leftskip 0pt plus 1fil\def\\{\par\smallskip}{\ifplaintex\large
\else\Large\fi\bf\thetitle}\par\medskip}   

\vglue 0pt plus 0.1fil 

{\parskip=0pt\leftskip 0pt plus 1fil\def\\{\par}{\sc\theauthors}
\par\medskip}

\vglue 0pt plus 0.1fil 

{\small\parskip=0pt\let\newline\\
{\leftskip 0pt plus 1fil\def\\{\par}{\sl\theaddress}\par}
\expandafter\ifx\theemail\relax    
\relax\else\vglue 5pt plus 0.02fil minus 2pt\def\\{\stdspace{\rm 
and}\stdspace} 
\cl{Email:\stdspace\tt\theemail}\fi
\ifx\theurl\relax                  
\relax\else\vglue 5pt plus 0.02fil minus 2pt\def\\{\stdspace{\rm 
and}\stdspace}
\cl{URL:\stdspace\tt\theurl}\fi\par}

\vglue 7pt plus 0.3fil minus 3pt

{\bf Abstract}
\vglue 5pt plus 0.1fil minus 2pt

\theabstract

\vglue 7pt plus 0.3fil minus 3pt

{\bf AMS Classification numbers}\quad Primary:\quad \theprimaryclass

Secondary:\quad \thesecondaryclass

\vglue 5pt plus 0.3fil minus 2pt

{\bf Keywords:}\quad \thekeywords

\vglue 10pt plus 0.5fil minus 5pt

{\small  Proposed: \theproposer\hfill Received: \receiveddate\nl
Seconded: \theseconders\hfill 
\ifx\reviseddate\relax                         
Accepted: \accepteddate                        
\else
Revised: \reviseddate                          
\fi}
\eject
}       

\let\maketitlepage\maketitlep
\let\maketitle\maketitlepage

\font\phead=cmsl9 scaled 950
\font=cmsl9 scaled 1050
\font\pnum=cmbx10 scaled 913
\font\lnum=cmbx10 
\font\pfoot=cmsl9 scaled 950
\font=cmsl9 scaled 1050
\ifplaintex
\headline{\vbox to 0pt{\vskip -4.5mm\line{\small\phead\ifnum
\count0=\startpage ISSN 1364-0380 (on line)
1465-3060 (printed) \hfill {\pnum\folio}\else\ifodd\count0\def\\{ }%
\ifx\theshorttitle\relax\thetitle\else\theshorttitle\fi\hfill{\pnum\folio}
\else\def\\{ and }{\pnum\folio}\hfill\ifx\theshortauthors\relax\theauthors
\else\theshortauthors\fi\fi\fi}\vss}}
\footline{\vbox to 0pt{\vglue 0mm\line{\small\pfoot\ifnum\count0=\startpage
\copyright\ \gtp\hfill\else
\gt, Volume \thevolumenumber\ (\thevolumeyear)\hfill\fi}\vss
}}
\else
\makeatletter
\def\@oddhead{{\small\lhead\ifnum\count0=\startpage ISSN 1364-0380 (on line)
1465-3060 (printed) \hfill {\lnum\number\count0}\else\ifodd\count0
\def\\{ }\ifx\theshorttitle\relax \thetitle \else\theshorttitle\fi\hfill
{\lnum\number\count0}\else\def\\{ and }{\lnum\number\count0}
\hfill\ifx\theshortauthors\relax 
\theauthors\else\theshortauthors\fi\fi\fi}}\def\@evenhead{@oddhead}
\def\@oddfoot{\small\lfoot\ifnum\count0=\startpage\copyright\ \gtp\hfill\else
\gt, Volume \thevolumenumber\ (\thevolumeyear)\hfill\fi}
\def\@evenfoot{@oddfoot}
\makeatother
\fi

\newwrite\gtoutfile
\long\gdef\makeheadfile{  
{\def\\{, }\def\s{ }
\immediate\openout\gtoutfile head.xxx
\immediate\write\gtoutfile{To: math@arxiv.org}
\immediate\write\gtoutfile{Subject: put OR rep NNNNN:pppp}
\immediate\write\gtoutfile{--text follows this line--}
\immediate\write\gtoutfile{Proxy-for: \ifx\theasciiauthors\relax
\theauthors\else\theasciiauthors\fi\s<\ifx\theasciiemail\relax\theemail\else\theasciiemail\fi>}
\immediate\write\gtoutfile{\noexpand\\}
\immediate\write\gtoutfile{Authors: \ifx\theasciiauthors\relax
\theauthors\else\theasciiauthors\fi}
{\def\\{ }\immediate\write\gtoutfile{Title: \ifx\theasciititle\relax
\thetitle\else\theasciititle\fi}}
\immediate\write\gtoutfile{Subj-class: GT or GR or SG or ...}
\immediate\write\gtoutfile{MSC-class: \theprimaryclass\ifx\thesecondaryclass\relax\else, \thesecondaryclass\fi}
\immediate\write\gtoutfile{Journal-ref: Geom. Topol. \thevolumenumber\s
(\thevolumeyear) \startpage-\finishpage}
\immediate\write\gtoutfile{Comments: Published in Geometry and Topology at}
\immediate\write\gtoutfile{    http://www.maths.warwick.ac.uk/gt/GTVol\thevolumenumber/paper\thepapernumber.abs.html}
\immediate\write\gtoutfile{\noexpand\\}
\immediate\write\gtoutfile{}
\ifx\theasciiabstract\relax
\immediate\write\gtoutfile{\theabstract}\else
\immediate\write\gtoutfile{\theasciiabstract}\fi
\immediate\write\gtoutfile{}
\immediate\write\gtoutfile{\noexpand\\}
\immediate\write\gtoutfile{}
\immediate\closeout\gtoutfile}}  

\def\maketitlepage{\maketitlep\makeheadfile}
\let\maketitle\maketitlepage

\lognumber{194}

\volumenumber{6}
\papernumber{20}
\volumeyear{2002}
\pagenumbers{563}{607}
\received{4 July 2001}
\accepted{28 October 2002}
\published{1 December  2002}
\proposed{Vaughan Jones}
\seconded{Robion Kirby, Joan Birman}

\usepackage{amsmath,amssymb,epsfig,longtable,rotating}  

\let\olditemize\itemize
\def\itemize{\olditemize\let\\\par}

\def\wh{\widehat}
\def\whP{\widehat\Phi_{l,E}}
\def\wbP{\overline\Phi_{l,E}}
\def\wb{\overline}

\def\Mod{\hbox{Mod}}

\def\spin{\hbox{\tiny{spin}}}
\def\ZZ{\mathbb{Z}}  
\def\QQ{\mathbb{Q}}
\def\RR{\mathbb{R}}
\def\CC{\mathbb{C}}
\def\NN{\mathbb{N}}
\def\go{\longrightarrow}
\def\im{\mapsto}
\def\inj{\hookrightarrow}
\def\A{{\cal{A}}}
\def\s{{\cal{S}}}
\def\M{{\cal{M}}}

\def\D{{\cal{D}}}

\def\E{{\cal{E}}}
\def\L{\Lambda}
\def\X{\chi}

\def\S{\mathfrak{S}}

\def\LL{{\cal{L}}}
\def\Ker{\hbox{Ker}\,}
\def\Im{\hbox{Im}\,}
\def\Hom{\hbox{Hom}}
\def\dim{\hbox{dim}}
\def\sdim{\hbox{sdim}}
\def\str{\hbox{str}}
\def\g{\mathfrak{g}}
\def\h{\mathfrak{h}}
\def\sll{\mathfrak{sl}}
\def\psl{\mathfrak{psl}}
\def\gl{\mathfrak{gl}}
\def\pgl{\mathfrak{pgl}}
\def\osp{\mathfrak{osp}}
\def\sp{\mathfrak{sp}}
\def\so{\mathfrak{so}}
\def\spin{\hbox{spin}}
\def\DD{{{\mathfrak D}_{2\,1}}}
\def\gg{{\mathfrak{g}}}
\def\ff{{\mathfrak{f}}}
\def\ee{{\mathfrak{e}}}
\def\ex{{\mathfrak{ex}}}
\def\Xd{\chi_{\DD}}
\def\K{\tilde{K}_3}

\def\labelb{\Phi_\DD(K) (x_1)=\sum_{\gamma_1,\gamma_2,\gamma_3}}
\newtheorem{prop}{Proposition}[section]
\newtheorem{lem}{Lemme}[section]
\newtheorem{theo}{Th{\'e}or{\`e}me}[section]
\newenvironment{dem}{\proof[D{\'e}monstration]}{\endproof}
\newenvironment{dem...}{\proof[D{\'e}monstration]}{}

\begin{document}
\title{Caract{\`e}res sur l'alg{\`e}bre de
      diagrammes trivalents $\L$}
\covertitle{Caract{\noexpand\`e}res sur l'alg{\noexpand\`e}bre de
      diagrammes trivalents $\Lambda$}
\asciititle{Caracteres sur l'algebre de
      diagrammes trivalents Lambda}

\author{Bertrand Patureau-Mirand}

\address{L.M.A.M. Universit{\'e} de Bretagne-Sud, Centre de 
Recherche\\Campus de Tohannic, BP 573, F-56017 Vannes, France}
\email{bertrand.patureau@univ-ubs.fr}
\url{http://www.univ-ubs.fr/lmam/patureau/}

\begin{abstract}
  The theory of Vassiliev invariants deals with many modules of
  diagrams on which the algebra $\L$ defined by Pierre Vogel in
  \cite{Vo2} acts. By specifying a quadratic simple Lie superalgebra,
  one obtains a character on $\L$. We show the coherence of these
  characters by building a map of graded algebras beetwen $\L$ and a
  quotient of a ring of polynomials in three variables; all the
  characters induced by simple Lie superalgebras factor through this
  map. In particular, we show that the characters for the Lie
  superalgebra $\ff(4)$ with dimension 40 and for $\sll_3$ are the
  same.
\smallskip

{\bf R\'esum\'e}

\smallskip
  De nombreux modules de diagrammes sont utilis{\'e}s dans la
  th{\'e}orie des invariants de Vassiliev. Pierre Vogel a d{\'e}finit
  dans \cite{Vo2} une alg{\`e}bre $\L$ qui agit sur ces espaces. Les
  superalg{\`e}bres de Lie simples quadratiques fournissent des
  caract{\`e}res sur $\L$. On montre leur coh{\'e}rence en
  construisant un morphisme d'alg{\`e}bre gradu{\'e}e, entre $\L$ et un
  quotient d'un anneau de polyn{\^o}me en trois variables, qui
  factorise tous ces caract{\`e}res. En particulier, on montre que le
  caract{\`e}re associ{\'e} {\'a} la superalg{\`e}bre de Lie $\ff(4)$ de
  dimension 40 co{\"\i}ncide avec celui associ{\'e} {\`a} $\sll_3$.
\end{abstract}

\asciiabstract{The theory of Vassiliev invariants deals with many
modules of diagrams on which the algebra Lambda defined by Pierre
Vogel acts.  By specifying a quadratic simple Lie superalgebra, one
obtains a character on Lambda. We show the coherence of these
characters by building a map of graded algebras beetwen Lambda and a
quotient of a ring of polynomials in three variables; all the
characters induced by simple Lie superalgebras factor through this
map. In particular, we show that the characters for the Lie
superalgebra f(4) with dimension 40 and for sl(3) are the same.
Resume: De nombreux modules de diagrammes sont utilises dans la
theorie des invariants de Vassiliev. Pierre Vogel a definit une
algebre Lambda qui agit sur ces espaces. Les superalgebres de Lie
simples quadratiques fournissent des caracteres sur Lambda.  On montre
leur coherence en construisant un morphisme d'algebre graduee, entre
Lambda et un quotient d'un anneau de polyneme en trois variables, qui
factorise tous ces caracteres. En particulier, on montre que le
caractere associe a la superalgebre de Lie f(4) de dimension 40
coincide avec celui associe a sl(3).}




\maketitlepage

\let\\\par

\section*{Introduction}
Cet article est tir{\'e} de mon travail en th{\`e}se. Il s'agit de la
d{\'e}monstration du th{\'e}or{\`e}me \ref{Pal} que j'ai annonc{\'e}
dans \cite{moi}.\\

L'int{\'e}grale de Kontsevich associe a un entrelacs son invariant de
Vassiliev--Kontsevich universel qui prend ses valeurs dans un espace vectoriel
engendr{\'e} par les diagrammes de cordes. Les espaces de diagrammes
trivalents (g{\'e}n{\'e}ralisant les diagrammes de cordes) ne sont connus que par une
pr{\'e}sentation. M{\^e}me les dimensions de ces espaces restent aujourd'hui inconnues.

D. Bar-Natan publie en 1995 un article dans lequel il utilise des alg{\`e}bres de
Lie quadratiques (munies de formes bilin{\'e}aires invariantes non d{\'e}g{\'e}n{\'e}r{\'e}es) pour
d{\'e}tecter des {\'e}l{\'e}ments des modules de diagrammes. Il construit des fonctions de
poids, qui sont des applications sur ces modules de diagrammes. Compos{\'e}es avec
l'invariant universel, elles donnent des invariants de type fini {\`a} valeurs dans
les tenseurs invariants d'une alg{\`e}bre de Lie.

La m{\^e}me ann{\'e}e, P. Vogel introduit des structures alg{\'e}briques suppl{\'e}men\-taires
sur ces modules de diagrammes et commence une {\'e}tude syst{\'e}matique de ces objets
et des fonctions de poids qui y sont d{\'e}finies. En particulier, il d{\'e}finit une
alg{\`e}bre $\Lambda$ qui agit sur plusieurs de ces modules. Les fonctions de poids
provenant de superalg{\`e}bres de Lie simples induisent des caract{\`e}res sur cette
alg{\`e}bre, {\`a} l'aide desquels il a {\'e}t{\'e} possible de montrer que les invariants de
type fini sont plus vastes que ceux correspondant aux invariants quantiques
classiques. La compr{\'e}hension de cette alg{\`e}bre est centrale pour la connaissance
globale des invariants de type fini.\\

Toute vari{\'e}t{\'e} de dimension trois peut {\^e}tre obtenue en faisant de la chirurgie le
long d'un entrelacs en bande. Cette description est utilis{\'e}e en 1995 par T.  Le,
H. et J. Murakami et T. Ohtsuki pour construire un invariant universel de type
fini pour les vari{\'e}t{\'e}s de dimension trois. Le logarithme de cet invariant prend
ces valeurs dans un espace isomorphe au compl{\'e}t{\'e} de l'alg{\`e}bre $\Lambda$. L'ann{\'e}e
suivante, T. Le et J.  Murakami donnent, en utilisant les travaux de Drinfield,
une version alg{\'e}brique de l'invariant universel de type fini.\\

Ce texte est organis{\'e} de la mani{\`e}re suivante:\\
Dans la premi{\`e}re partie, j'introduit les modules de diagrammes et je rappelle
leur lien avec l'int{\'e}grale de Kontsevich.\\ 
Dans la deuxi{\`e}me partie, j'introduit les fonctions de poids associ{\'e}es aux
superalg{\`e}bres de Lie et j'{\'e}nonce le th{\'e}or{\`e}me \ref{Pal} sur la coh{\'e}rence des
caract{\`e}res. Sa d{\'e}monstration repose sur la construction, dans la troisi{\`e}me
partie, d'un morphisme d'alg{\`e}bre gradu{\'e}e entre $\Lambda$ et un quotient d'un
anneau de polyn{\^o}me {\`a} trois variables qui factorise tous les caract{\`e}res induits
par les superalg{\`e}bres de Lie simples.\\
Les cas des superalg{\`e}bres $\gg(3)$ et $\ff(4)$ est trait{\'e} s{\'e}par{\'e}ment dans la
quatri{\`e}me partie.

\rk{Remerciements}
Je remercie P. Vogel qui a suivi ce travail durant ma th{\`e}se de
doctorat et le r{\'e}f{\'e}r{\'e} {\`a} qui est d{\^u} cet index des
notations.

\rk{\large \bf Index des notations}

\smallskip\small
\begin{longtable}{p{2.5cm}p{1cm}p{7.5cm}}
$\Gamma$\dotfill&&d{\'e}signe une vari{\'e}t{\'e} de dimension $1$ {\`a} bord\\
$X$\dotfill&&d{\'e}signe un ensemble fini \\
$[n]$\dotfill&&ensemble des entiers de $1$ {\`a} $n$\\
$\S(X)$, $\S_n$\dotfill&&groupe des permutations de $X$ (resp.\ de $[n]$)\\
$\A(\Gamma,X)$\dotfill&(\ref{AGX})& espace des $(\Gamma,X)$--diagrammes\\
$F(X),\,F_n$\dotfill&(\ref{AGX})& espace des $(\emptyset,X)$--diagrammes connexes ($F_n=F([n])$)\\
$\L$\dotfill&(\ref{Lambda})&alg{\`e}bre des diagrammes ``{\`a} $3$ jambes''\\
$\L_0,\,t,\,x_n$\dotfill&(\ref{Lambda})&$\L_0$ est la sous alg{\`e}bre de $\L$ engendr{\'e}e par les {\'e}l{\'e}ments $t$, $x_n$\\
$\X_\square$\dotfill&&les caract{\`e}res (morphismes d'alg{\`e}bres gradu{\'e}s) sur $\L$ sont not{\'e}s par la lettre $\X$\\
$\D$\dotfill&(\ref{catD})&cat{\'e}gorie des $(\emptyset,X)$--diagrammes\\
$\D_\Gamma$\dotfill&(\ref{catD})&cat{\'e}gorie des $((S^1)^{\amalg n},X)$--diagrammes\\
$\D_b$\dotfill&(\ref{catD})&cat{\'e}gorie dans laquelle les morphismes sont les diagrammes connexes relativement au but\\
$\D\stackrel{\Phi_L}{\go} \Mod_L$\dotfill&(\ref{PhiL})&foncteur mono{\"\i}dal lin{\'e}aire associ{\'e} {\`a} $L$\\
$S=\QQ[t,u,v]$\dotfill&(\ref{annS})&anneau gradu{\'e} contenant les polyn{\^o}mes $P_{\sll},P_{\osp},P_{\DD},P_{\sll_2},P_{\ex},Q_{\ex}$\\
$\D_\osp,\D_\sll,\D_\gl$, $\Phi_\osp,\Phi_\sll,\Phi_\gl $\dotfill&(\ref{DM})&cat{\'e}gories quotientes de $\D$ et leur foncteur quotient\\
$\D_{\gl_0}$, $\Phi_{\gl_0}$\dotfill&(\ref{Dgl0})&variante de $\D_\gl$ correspondant a une superdimension nulle\\
$\M,\, \Phi_\M$\dotfill&(\ref{DM})&cat{\'e}gorie des $X$--surfaces marqu{\'e}es et le foncteur d'{\'e}pai\-cis\-se\-ment $\D\go\M$\\
$\s$\dotfill&(\ref{anns})&alg{\`e}bre du mono{\"\i}de des surfaces compactes $\s\simeq\QQ[\delta,\alpha,\beta]_{/(\alpha\beta-\alpha^3)}$\\
$M(X)$, $M[n]$, $M_c(X)$\dotfill&(\ref{MX})&espace des $X$-- (resp.\ $[n]$--) surfaces marqu{\'e}es (resp.\ connexes)\\
$R,\,\sigma_2,\,\sigma_3$\dotfill&(\ref{D21})&$R=\QQ[a,b,c]_{/(a+b+c)}$ contient les {\'e}l{\'e}ments $\sigma_2=ab+bc+ca$ et $\sigma_3=abc$\\
$\DD$\dotfill&(\ref{D21})&$R$--superalg{\`e}bre de Lie g{\'e}n{\'e}rique pour les superalg{\`e}bres ${\mathfrak D}_{2\,1,\alpha}$\\
$\wb\DD$\dotfill&(\ref{D21})&$\wb\DD=\DD\otimes_Rk$ o{\`u} $k$ est le corps des fractions de $R$\\
$N=\NN\amalg N_3\amalg N_6$\dotfill&(\ref{NN})&ensembles de multi-entiers munis d'un bon ordre\\
$W_\delta$\dotfill&(\ref{NN})&diagrammes, form{\'e}s de g{\'e}n{\'e}ralisations des roues, param{\'e}tr{\'e}s par $N$\\
$f_K$\dotfill&(\ref{NN})&application obtenue par recollement du diagramme $K$\\
$R_\delta,\,\bar R_\delta$\dotfill&(\ref{NN})&filtration de $F_0$ index{\'e}e par $N$ et gradu{\'e} associ{\'e}\\
$\wh\A(X_1,X_2)$, $\wb\A(X_1,X_2))$\dotfill&(\ref{Dbic})&espace des diagrammes bicolores (resp.\ modulo la relation ($\wb{IHX}$))\\
$\wh\D$, $\wb\D$\dotfill&(\ref{Dbic})&cat{\'e}gorie des diagrammes bicolores (resp.\ modulo la relation ($\wb{IHX}$))\\
\end{longtable}\normalsize

\section{Les diagrammes trivalents}
\subsection{Modules de diagrammes}
Dans tout ce qui suit, $\Gamma$ sera une vari{\'e}t{\'e} compacte de dimension un {\`a} bord
et $X$ d{\'e}note un ensemble fini. Un $(\Gamma,X)$--diagramme est un graphe fini $K$,
dont tous les sommets sont trivalents ou monovalents, muni des donn{\'e}es
suivantes:
\begin{enumerate}
\item Un isomorphisme d'un sous-graphe de $K$ vers $\Gamma \amalg X$ envoyant
  l'ensemble des sommets monovalents de $K$ sur $\partial\Gamma\cup X$.
\item Pour chaque sommet trivalent $x$ de $K$, un ordre cyclique sur
  l'ensemble des trois ar{\^e}tes orient{\'e}es arrivant en $x$.
\end{enumerate}                                
On peut repr{\'e}senter un $(\Gamma,X)$--diagramme par un graphe trivalent immerg{\'e}
dans le plan de mani{\`e}re {\`a} ce que l'ordre cyclique en chaque sommet soit donn{\'e} par
l'orientation du plan. On repr{\'e}sentera d'un trait plus {\'e}pais les ar{\^e}tes
appartenant {\`a} $\Gamma$.\\
Soit $E$ le quotient du $\QQ$--espace vectoriel librement engendr{\'e} par les diagrammes
trivalents par les relations suivantes:
\begin{enumerate}
\item Si deux diagrammes ne diff{\`e}rent que par l'ordre cyclique de l'un
  de leurs sommets, leur somme est nulle (relation dite (AS) pour
antisym{\'e}trie):
  $$\begin{array}{cccc}\put(-10,-10) {\epsfbox{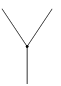}}&+&
    \put(-10,-10) {\epsfbox{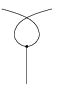}}&\equiv0\end{array}$$
\item La relation (IHX) fait intervenir trois diagrammes qui ne
  diff{\`e}rent qu'au voisinage d'une ar{\^e}te:
  $$\begin{array}{ccccc}\put(-8,-5) {\epsfbox{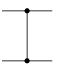}}&\equiv&
    \put(-8,-5) {\epsfbox{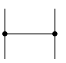}} &
    -&\put(-8,-5) {\epsfbox{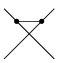}}\end{array}$$
\item La relation (STU) qui est une variation de la relation (IHX) au voisinage de
  $\Gamma$:                               
$$\begin{array}{ccccc}\put(-8,-5) {\epsfbox{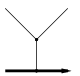}}&\equiv&
\put(-8,-5) {\epsfbox{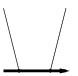}} & -&\put(-8,-5) {\epsfbox{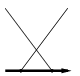}}\end{array}$$
\end{enumerate}

On d{\'e}signe par $\A(\Gamma,X)$ \label{AGX} le sous-espace de $E$ engendr{\'e} par les
$(\Gamma,X)$--diagram\-mes et par $F(X)$ le sous-espace de $E$ engendr{\'e} par les
$(\emptyset,X)$--diagrammes connexes ayant au moins un sommet trivalent.\\
Enfin, on note $[n]$ l'ensemble $\{1,2,\ldots,n\}$ et $F_n$ pour $F([n])$.

On d{\'e}finit le degr{\'e} d'un diagramme $K\in E$ par $a-s$ o{\`u} $a$ et $s$
sont les nombres d'ar{\^e}tes et de sommets trivalents de $K$. Ainsi, ces modules
sont munis d'une graduation. On note $d$ le diagramme de
$\A(\emptyset,\emptyset)$ form{\'e} d'un seul cercle et on conviendra que son degr{\'e}
est nul.

Toute bijection entre des ensembles finis $X$ et $Y$ induit une bijection entre
$\A(\Gamma,X)$ et $\A(\Gamma,Y)$. En particulier le groupe sym{\'e}trique $\S_n$
op{\`e}re sur $F_n$.

On dira qu'un diagramme $K$ a $n$ boucles si la dimension de son premier groupe
d'homologie est $n$ (i.e. $\dim(H_1(K))=n$). Si $n$ est un entier strictement
positif, le degr{\'e} d'un diagramme de $F_n$ est {\'e}gal {\`a} son nombre de boucles plus
$n-1$.

On a ici repris les notations de \cite{Vo2} {\`a} l'exception des coefficients qui sont
ici rationnels, de la d{\'e}finition du degr{\'e} et des d{\'e}finitions de $F_n$ qui
n'entra{\^\i}nent des modifications que pour $F_0$ et $F_2$ qui sont ici pris nuls en
degr{\'e}s respectifs z{\'e}ro et un.

\subsection{L'int{\'e}grale de Kontsevich}
Il est connu que le module gradu{\'e} $\A=\A(S^1,\emptyset)$ peut {\^e}tre
muni d'une structure d'alg{\`e}bre de Hopf gradu{\'e}e, commutative et cocommutative.\\
L'alg{\`e}bre $\A$ est donc l'alg{\`e}bre sym{\'e}trique du sous-module gradu{\'e} $\cal{P}$
form{\'e} par ses {\'e}l{\'e}ments primitifs et ce module est reli{\'e} aux modules $F_n$ par
l'isomorphisme:
$${\cal P}\simeq\bigoplus_{n=2}^{+\infty} H^0(F_n,\S_n)$$  
Si on note $\Theta$ le $(S^1,\emptyset)$--diagramme repr{\'e}sent{\'e} par le cercle et
un de ses diam{\`e}tres, alors l'int{\'e}grale de Kontsevich associe {\`a} chaque n{\oe}ud
orient{\'e} son invariant de Vassiliev universel {\`a} valeurs dans l'alg{\`e}bre $\hat\A_r$
qui est la compl{\'e}t{\'e}e pour la graduation du quotient $\A_{/(\Theta)}$.\\  

\subsection{L'alg{\`e}bre $\L$}
Dans cette section est introduite l'alg{\`e}bre gradu{\'e}e de diagramme $\L$ qui agit
de mani{\`e}re naturelle sur les modules $F(X)$.\\ 
$\L$ \label{Lambda} est d{\'e}finie comme le sous-espace vectoriel
form{\'e} des {\'e}l{\'e}ments de $F_3$ totalement antisym{\'e}triques
sous l'action du groupe sym{\'e}trique.\\
$\L$ est naturellement munie d'une structure d'alg{\`e}bre commutative et agit sur
chaque module $F(X)$. Si $u$ appartient {\`a} $\L$ et $K\in F(X)$
est un diagramme, un diagramme repr{\'e}sentant $u.K$ est obtenu en ins{\'e}rant $u$ au
niveau d'un sommet trivalent de $K$. Afin de rendre cette action coh{\'e}rente avec
les graduations des modules, on convient de d{\'e}finir le degr{\'e} d'un {\'e}l{\'e}ment de
$\L$ comme son degr{\'e} dans $F_3$ moins deux (de sorte que l'unit{\'e} de $\L$ est de
degr{\'e} nul).\\
On a la description suivante de $F_n$ pour $n$ petit (cf \cite{Vo2}): $F_0$ est un
$\L$--module libre de rang un engendr{\'e} par l'unique diagramme (aux relations AS
pr{\`e}s) de degr{\'e} un; $F_1$ est nul; $F_2$ est un $\L$--module libre de rang un
engendr{\'e} par l'unique diagramme (aux relations AS pr{\`e}s) de degr{\'e} deux. De plus
on ne conna{\^\i}t pas d'exemple d'{\'e}l{\'e}ment de $F_3$ qui ne soit pas dans
$\L$. D'autre part, $\L$ est engendr{\'e}e en degr{\'e} $1$ par l'{\'e}l{\'e}ment $t$ ci dessous
et on peut construire la famille $x_n$ d'{\'e}l{\'e}ments de $\L$ qui engendrent avec
$t$ une sous-alg{\`e}bre de $\L$ not{\'e}e $\L_0$ .
$$\begin{array}{cccccc}t=&\put(-8,-5) {\epsfbox{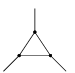}}&=\frac12&
\put(-8,-5) {\epsfbox{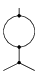}} & \qquad x_n=&\put(-8,-10) {\epsfbox{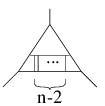}}
\end{array}$$
Il n'y a pas de contre-exemple {\`a} la conjecture suivante: $\L=\L_0$.
C'est la motivation principale de cet article.

\section{Fonction de poids associ{\'e}e {\`a} une superalg{\`e}bre de Lie quadratique}
\subsection{Le foncteur $\Phi_L$}
Dans cette partie, nous introduisons les fonctions de poids g{\'e}n{\'e}ralis{\'e}es
associ{\'e}es aux superalg{\`e}bres de Lie quadratique. Ces applications sont
aujourd'hui les seules mani{\`e}res connues de d{\'e}tecter la non nullit{\'e} des {\'e}l{\'e}ments
des modules de diagrammes.

Soit $L$ une superalg{\`e}bre de Lie sur un corps $k$ de caract{\'e}ristique nulle munie
d'un {\'e}l{\'e}ment de Casimir non d{\'e}g{\'e}n{\'e}r{\'e} $\Omega\in L\otimes L$ de degr{\'e} pair.  Le
Casimir fournit un isomorphisme de $L$--module entre $L$ et son dual et la forme
bilin{\'e}aire sur $L$ supersym{\'e}trique invariante associ{\'e}e sera not{\'e}e $<.,.>$.  On
construit une cat{\'e}gorie $\D$ de diagrammes et un foncteur, not{\'e}
$\Phi_{L,\Omega}$, de $\D$ vers la cat{\'e}gorie $\Mod_L$ des repr{\'e}sentations de
$L$.\\
\medskip
Soit $\D$ \label{catD} la cat{\'e}gorie $\QQ$--lin{\'e}aire mono{\"\i}dale d{\'e}finie par:
$$\hbox{Obj}(\D)=\left\{[n],n\in\NN\right\}$$
$$\D([p],[q])=\A(\emptyset,[p]\amalg[q])$$ 
La composition d'un diagramme de $[p]$ vers $[q]$ avec un diagramme de $[q]$
vers $[r]$ est donn{\'e}e par la r{\'e}union au dessus de $[q]$ des deux diagrammes (on
les recolle en identifiant les sommets monovalents de m{\^e}me index des deux
ensembles $[q]$). On {\'e}tend cette d{\'e}finition par lin{\'e}arit{\'e} {\`a} des combinaisons
lin{\'e}aires de diagrammes.\\
Le produit tensoriel $[p]\otimes [q]$ vaut $[p+q]$ et celui de deux diagrammes
est donn{\'e} par l'image de leur r{\'e}union disjointe par l'isomorphisme de $[p]\amalg
[q]\simeq [p+q]$ obtenu en augmentant de $p$ chaque {\'e}l{\'e}ment de $[q]$.\\ 
\medskip
\label{PhiL}
\begin{prop} {\rm(cf \cite{Vo2})}\\
Il existe un unique foncteur $\QQ$--lin{\'e}aire mono{\"\i}dal $\Phi_{L,\Omega}$
de la cat{\'e}gorie $\D$ vers la cat{\'e}gorie $\Mod_L$ envoyant $[n]$ sur
$L^{\otimes n}$ et les diagrammes suivants:
$${\epsfbox{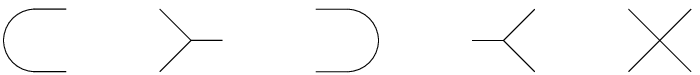}}$$
vers respectivement:
\begin{enumerate}
\item Le Casimir $\Omega\in L^{\otimes 2}\simeq \Mod_L(\QQ,L^{\otimes 2})$
\item Le crochet de Lie de $L^{\otimes 2}$ vers $L$
\item Le produit scalaire associ{\'e} au Casimir de $L^{\otimes 2}$ vers
  $\QQ$
\item Le dual du crochet de Lie de $L$ vers $L^{\otimes 2}$
\item L'op{\'e}rateur de sym{\'e}trie: $\begin{array}[t]{ccl} L^{\otimes
      2}&\go&L^{\otimes 2}\\ x\otimes
    y&\im&(-1)^{\hbox{degr{\'e}}(x)\hbox{degr{\'e}}(y)}y\otimes x
  \end{array}$
\end{enumerate}
D'autre part, si $L$ est simple, il existe un caract{\`e}re gradu{\'e}:
$\X_L\co \L\go\QQ[x]$ v{\'e}rifiant: 
$$\forall u\in\L,\, \forall K\in F([p]\amalg[q])\subset\D([p],[q]),\,
\Phi_{L,\Omega}(uK)=\X_L(u)_{|x=1}\Phi_{L,\Omega}(K)$$
\end{prop}
\paragraph{Remarque} On peut d{\'e}finir la notion de ``pseudo-alg{\`e}bre de Lie''
comme une cat{\'e}gorie $\LL$ $\QQ$--lin{\'e}aire mono{\"\i}dale munie d'un foncteur $\Phi_\LL$
de $\D$ vers $\LL$. Ceci signifie que $\LL$ poss{\`e}de un objet particulier not{\'e} $L$,
un op{\'e}rateur de sym{\'e}trie (endomorphisme de $L^{\otimes 2}$) induisant une
repr{\'e}sentation de groupe sym{\'e}trique $\S_n$ dans les endomorphismes de
$L^{\otimes n}$, un op{\'e}rateur de Casimir sym{\'e}\-tri\-que ayant pour adjoint un morphisme
$<.,.>\co L^{\otimes 2}\go 1_\otimes$ (dans le sens o{\`u} $L\simeq L\otimes
1_\otimes\stackrel{Id_L\otimes\Omega}{\go}L^{\otimes 3}\stackrel{<.,.>\otimes
  Id_L}{\go}1_\otimes \otimes L\simeq L\equiv Id_L$) et un op{\'e}rateur ``crochet de
Lie'' $[.,.] \co  L^{\otimes 2}\go L$ antisym{\'e}trique et
v{\'e}rifiant l'identit{\'e} de Jacobi.\\
La proposition dit alors que la cat{\'e}gorie des repr{\'e}sentations d'une superalg{\`e}bre
de Lie quadratique a naturellement une structure de pseudo-alg{\`e}bre de Lie.\\

Si $f$ est une bijection de $[p]\amalg [q]$ vers $[r]\amalg [s]$, alors $f$
induit un isomorphisme de $\Hom_L(L^{\otimes p},L^{\otimes q})$ vers
$\Hom_L(L^{\otimes r},L^{\otimes s})$ (par l'autodualit{\'e} de $L$) mais aussi un
isomorphisme de $\A(\emptyset,[p]\amalg[q])$ vers $\A(\emptyset,[r]\amalg[s])$
et il est facile de voir que $f_*\circ\Phi_L=\Phi_L\circ f_*$. En particulier,
si $p+q=n$, on identifie souvent $F_n$ {\`a} un sous-module de
$\A(\emptyset,[p]\amalg[q])$. \\

On peut d{\'e}finir des variantes du foncteur $\Phi_L$: 
\begin{itemize}
\item On d{\'e}finit $\D_\Gamma$ comme la cat{\'e}gorie ayant les m{\^e}mes objets que $\D$
  mais dont les morphismes sont donn{\'e}s par
  $$\D_\Gamma([p],[q])=\bigoplus_{n\in\NN}\A((S^1)^{\amalg n},[p]\amalg[q])$$ 
  Si d'autre part $E$ est un $L$--module, on peut encore montrer qu'il existe un
  unique foncteur $\Phi_{L,E}\co \D_\Gamma\go\Mod_L$ prolongeant $\Phi_L$ et
  envoyant
  $${\epsfbox{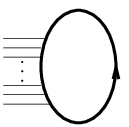}}$$
  sur $x_1\otimes x_2 \otimes \ldots \otimes
  x_n\mapsto \str_E(x_1x_2\ldots x_n)$. La notation $\str_E$ d{\'e}signe la
  supertrace sur le supermodule $E$ et, par la suite, on notera $\sdim(E)$ la
  superdimension d'un module $E$ (qui est {\'e}gale {\`a} la dimension de la partie
  paire de $E$ moins la dimension de sa partie impaire). 
\item Si $L$ est une superalg{\`e}bre de lie sur une $\QQ$--alg{\`e}bre $R$ munie d'un {\'e}l{\'e}ment
  de Casimir $\Omega\in L\otimes_R L$ et si l'on consid{\`e}re $\D_b$ la
  sous-cat{\'e}gorie de $\D$ ayant les m{\^e}mes objets et dont les morphismes sont
  engendr{\'e}s par les $(\emptyset,[p]\amalg[q])$--diagrammes dont toutes les
  composantes connexes rencontrent $[q]$, alors on peut d{\'e}finir de m{\^e}me un
  foncteur $R$--lin{\'e}aire: $\Phi_L\co \D_b\go\Mod_L$ qui co{\"\i}ncide avec la restriction
  {\`a} $\D_b$ du foncteur $\Phi_L$ de la proposition si $R$ et $L$ satisfont aux
  hypoth{\`e}ses.
\end{itemize}
Enfin on se servira du lemme:
\begin{lem}\label{Phi_b}
  Soit $\psi\co R\go R'$ un morphisme entre deux $\QQ$--alg{\`e}bres commutatives, et
  $f\co \h\go \g$ un morphisme entre la $R$--superalg{\`e}bre de Lie $\h$ et la
  $R'$--superalg{\`e}bre de Lie $\g$; par $\psi$, $\g$ est naturellement munie d'une
  structure de $R$--module et on suppose que $f$ est un morphisme de $R$--alg{\`e}bre
  de Lie. Si $\h$ poss{\`e}de un {\'e}l{\'e}ment de Casimir invariant $\Omega\in
  \h\otimes\h$ qui est envoy{\'e} par $f$ sur $f_*(\Omega)\in \g\otimes\g$ {\'e}l{\'e}ment
  de Casimir $\g$--invariant de $\g$, alors l'application $f_*$ induite par $f$
  et $\psi$ entre la $R$--alg{\`e}bre tensorielle de $\h$ et la $R'$--alg{\`e}bre
  tensorielle de $\g$ v{\'e}rifie pour tout $K\in\D_b([p],[q])$,
  $$f_*\circ\Phi_{\h,\Omega}(K)=\Phi_{\g,f_*\Omega}(K)\circ f_*$$
\end{lem}
Ce lemme est une cons{\'e}quence du fait que $\D_b$ est engendr{\'e}e comme cat{\'e}gorie
$\ZZ$--lin{\'e}aire mono{\"\i}dale par les morphismes ``crochet'', ``Casimir'' et
``sym{\'e}trie'' qui v{\'e}rifient tous le lemme par hypoth{\`e}se.

\subsection{Propri{\'e}t{\'e}s communes {\`a} tous les foncteurs $\Phi_L$}
On conna{\^\i}t la liste compl{\`e}te des superalg{\`e}bres de Lie simples complexes
quadratiques (cf \cite{Kac} et \cite{Vo2}). Elle est form{\'e}e des superalg{\`e}bres
$\sll(V)$ o{\`u} $V$ est un superespace de superdimension non nulle, des
superalg{\`e}bres $\psl(V)$ si $V$ est de superdimension nulle, des superalg{\`e}bres
$\osp(V)$ lorsque $V$ est muni d'une forme bilin{\'e}aire supersym{\'e}trique non
d{\'e}g{\'e}n{\'e}r{\'e}e, les superalg{\`e}bres $\DD_{,\delta}$ o{\`u} $\delta$ est un param{\`e}tre
complexe diff{\'e}rent de $0$ et de $1$, les cinq alg{\`e}bres exceptionnelles, les deux
superalg{\`e}bres $\gg(3)$ et $\ff(4)$ et enfin les superalg{\`e}bres
hamiltoniennes. Nous exclurons ces derni{\`e}res qui induisent sur $\L$ le caract{\`e}re
trivial (cf \cite{Vo2}, derni{\`e}re remarque de la partie 6).\\

L'action du Casimir $\Omega$ de $L$ sur $L^{\otimes n}$ s'exprime comme l'image
par $\Phi_L$ d'un diagramme de $\A(\emptyset,[n]\amalg[n])$. Sous l'action du
Casimir, $L^{\otimes n}$ se scinde en espaces caract{\'e}ristiques. En particulier,
il est facile de voir que $\Omega$ agit par $2t$ sur $L$ (o{\`u} l'on note encore
$t$ l'{\'e}l{\'e}ment $\X_L(t)_{|x=1}$).\\
Sur $L^{\otimes 2}$, $\Omega$ agit comme $4t-2\Psi$ o{\`u}: 
$$\begin{picture}(0,0)%
\epsfig{file=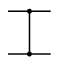}%
\end{picture}%
\setlength{\unitlength}{0.000625in}%
\begingroup\makeatletter\ifx\SetFigFont\undefined%
\gdef\SetFigFont#1#2#3#4#5{%
  \reset@font\fontsize{#1}{#2pt}%
  \fontfamily{#3}\fontseries{#4}\fontshape{#5}%
  \selectfont}%
\fi\endgroup%
\begin{picture}(1963,340)(108,-105)
\put(2026, 14){\makebox(0,0)[lb]{\smash{\SetFigFont{14}{16.8}{\rmdefault}{\mddefault}{\updefault}$)$}}}
\put(1801, 14){\makebox(0,0)[rb]{\smash{\SetFigFont{14}{16.8}{\rmdefault}{\mddefault}{\updefault}$\Psi=\Phi_L($}}}
\end{picture}
$$
On peut toujours d{\'e}composer $L^{\otimes 2}$ de la mani{\`e}re suivante:
$$L^{\otimes 2}=\L^2L\oplus S^2L$$ De plus si $t\neq0$, $\L^2L$ se d{\'e}compose en
la somme directe $L\oplus X_2$, $\Psi$ agit par $t$ sur $L$ et par $0$ sur
$X_2$. Si $L$ est l'une des superalg{\`e}bre de Lie de la famille ${\mathfrak
D}_{2\,1,\delta}$ alors $t=0$ et $\L^2L$ est une extension g{\'e}n{\'e}ralement non
scind{\'e}e de $X_2$ (noyau du crochet de Lie) par $L$.\\

Enfin, dans tous les cas, $S^2L$ s'{\'e}crit $\QQ\oplus E$ (le module trivial {\'e}tant
engendr{\'e} par le Casimir) et le module $E$ se d{\'e}compose en trois espaces propres,
les valeurs propres associ{\'e}es ${\alpha,\beta,\gamma}$ ont pour somme
$t$. Certains des modules cit{\'e}s ci-dessus peuvent {\'e}ventuellement {\^e}tre nuls mais
il existe toujours des {\'e}l{\'e}ments
$u(L)=-\frac{\alpha\beta+\beta\gamma+\gamma\alpha}{2}$ et
$v(L)=\frac{\alpha\beta\gamma}{2}$ tels que sur $E$ on ait la relation
\begin{equation}\label{psi3}
  \Psi^3=t\Psi^2+2u\Psi+2v.
\end{equation}
Le triplet $(\alpha,\beta,\gamma)\in\CC^3$ sera dit admissible pour $(L,\Omega)$
si l'{\'e}quation (1) est v{\'e}rifi{\'e}e. Si un triplet est admissible pour $L$, tout
triplet obtenu par permutation des trois valeurs $(\alpha,\beta,\gamma)$ est
bien s{\^u}r aussi admissible. De m{\^e}me, tout triplet non nul proportionnel {\`a} ce
triplet sera admissible pour un autre choix du Casimir.\\ 
Si $L$ n'est pas isomorphe {\`a} l'alg{\`e}bre de Lie $\sll_2$, et {\`a} permutation des
trois valeurs pr{\`e}s, l'ensemble des triplets admissibles pour $L$ sont sur une
unique droite de l'espace affine $\CC[\alpha,\beta,\gamma]$. Il leur correspond
donc un unique point de ${\mathbb{P}}(\CC[\alpha,\beta,\gamma])$. Pour $\sll_2$,
toujours {\`a} permutation pr{\`e}s des trois valeurs, l'ensemble des triplets
admissibles forme une droite du plan projectif complexe
${\mathbb{P}}(\CC[\alpha,\beta,\gamma])$.\\
Des superalg{\`e}bres distinctes peuvent avoir les m{\^e}mes triplets admissibles; c'est
le cas pour $\sll_3$ et $\ff(4)$, et de m{\^e}me, les triplets admissibles pour
$\gg(3)$ le sont aussi pour $\sll_2$.\\
On peut remarquer que l'ensemble de tous les triplets admissibles se trouve sur
cinq r{\'e}unions de droites de ${\mathbb{P}}(\CC[\alpha,\beta,\gamma])$:
\begin{itemize}
\item L'alg{\`e}bre $\sll_2$ est un cas particulier: $E$ est alors simple et $\Psi$
y vaut $-t$. Tous les triplets de la forme $(-t,\beta,2t-\beta)$ sont
admissibles pour $\sll_2$. On a toujours, {\`a} permutation des racines de $\Psi$
pr{\`e}s, $t+\alpha=0$ et le polyn{\^o}me
$P_{\sll_2}=\frac12(t+\alpha)(t+\beta)(t+\gamma)=v-ut+t^3$ est toujours nul.
\item Si $L$ est l'alg{\`e}bre $\sll(V)$ o{\`u} $V$ est un supermodule de superdimension
$\delta$, alors les triplets de la forme $(2,-2,\delta)$ sont admissibles pour
pour $L$. A permutation des racines de $\Psi$ pr{\`e}s, on a donc toujours
$\alpha+\beta=0$, donc le polyn{\^o}me d{\'e}fini par $P_{\sll}=v+ut$ est nul.
\item Si $L$ est l'alg{\`e}bre $\osp(V)$ o{\`u} $V$ est un supermodule muni d'une forme
  bilin{\'e}aire supersym{\'e}trique non d{\'e}g{\'e}n{\'e}r{\'e}e de superdimension $\delta$, alors les
  triplets de la forme $(4,-2,\delta-4)$ sont admissibles pour $L$. A
  permutation des racines de $\Psi$ pr{\`e}s, on a $\alpha+2\beta=0$, donc le
  polyn{\^o}me d{\'e}fini par $P_{\osp}=27v^2+18vut+2vt^3-8u^3+8u^2t^2$ est nul.
\item Si $L$ est une des superalg{\`e}bres de Lie ${\mathfrak D}_{2\,1,\delta}$ o{\`u}
  $\delta$ est un nombre complexe quelconque, alors le triplet
$(1,\delta,-1-\delta)$ est admissible pour $L$. La somme $t$ des racines de
  $\Psi$ est nulle et on pose $P_{\DD}=t$.
\item Si $L$ est une alg{\`e}bre de Lie exceptionnelle, $\Psi$ a deux valeurs
  propres de somme $\frac{t}{3}$. Les triplets admissibles pour $L$ sont donc de
  la forme $(\frac{2t}{3},\alpha(L),\frac{t}{3}-\alpha(L))$ et le polyn{\^o}me
  $P_{\ex}=27v+18ut+2t^3$ est nul.
\end{itemize}

Dans la suite, on note $S=\QQ[t,u,v]$ \label{annS}l'anneau des polyn{\^o}mes {\`a} trois
ind{\'e}termi\-n{\'e}es de degr{\'e} respectif $1$, $2$ et $3$. Pour chaque superalg{\`e}bre de
Lie quadratique munie d'un triplet admissible $L$, on pose $I_L=Ker(f_L)$ o{\`u}
$f_L$ est le morphisme d'alg{\`e}bre gradu{\'e}e de $S$ dans $\QQ[x]$ qui envoie les
ind{\'e}termin{\'e}es $t$, $u$ et $v$ sur leurs valeurs respectives $t(L)x$, $u(L)x^2$
et $v(L)x^3$ dans $\QQ[x]$. Nous allons montrer le th{\'e}or{\`e}me:
\begin{theo} \label{Pal}
Soit $I=\bigcap_L I_L$ alors il existe un unique caract{\`e}re gradu{\'e}\\ 
$\X\co \L\go S/I$ tel que pour toute superalg{\`e}bre de Lie $L$ de la liste ci-dessus,
on ait:
$$\X_L=f_L\circ\X$$
\end{theo}

\paragraph{\bf Remarque} 
L'id{\'e}al $I$ du th{\'e}or{\`e}me est la somme des deux id{\'e}aux principaux engendr{\'e}s par les
polyn{\^o}mes $P$ et $Q$ suivants:
$$P=P_{\sll}P_{\osp}P_{\DD}P_{\sll_2}P_{\ex}$$
$$Q=P_{\sll}P_{\osp}P_{\DD}P_{\sll_2}Q_{\ex}$$
o{\`u} $Q_{\ex}$ est le polyn{\^o}me de degr{\'e} dix d{\'e}fini section \ref{theo d'exis}.
D'autre part les calculs de \cite{Kn1} montrent l'existence d'un morphisme d'alg{\`e}bre
surjectif:
$$S_0=\QQ[t]\oplus P_{\sll_2}.S\go\L_0$$ 
dont la compos{\'e}e avec $\X$ est le morphisme quotient.\\
Les calculs par ordinateur de \cite{Kn1} montrent en outre que l'application
$S_0\go\L$ est un isomorphisme en degr{\'e} inf{\'e}rieur ou {\'e}gal {\`a} dix.

\rk{Remarques sur $\L$}
 
La nullit{\'e} d'une combinaison lin{\'e}aire de diagrammes est tr{\`e}s difficile {\`a} d{\'e}tecter
car la taille des pr{\'e}sentations des modules de diagrammes cro{\^\i}t tr{\`e}s vite avec
le degr{\'e}. Ainsi, lorsque P. Vogel a construit les {\'e}l{\'e}ments $x_n$ qui forment la
sous-alg{\`e}bre $\L_0$ de $\L$, il a montr{\'e} que ces {\'e}l{\'e}ments pris pour $n$ impair
suffisent {\`a} engendrer $\L_0$ et il n'y avait pas {\`a} priori de raison de
supposer l'existence d'autres relations entre les $x_n$. J.A. Kneissler a
ensuite montr{\'e} l'existence de relations suppl{\'e}mentaires permettant de construire
le morphisme d'alg{\`e}bre $S_0\go\L_0$. Cela donnait {\`a} penser que $\L$ aurait pu
{\^e}tre isomorphe {\`a} l'anneau $S_0$. Ceci aussi s'est r{\'e}v{\`e}l{\'e} inexact:\\ 
Les calculs r{\'e}cents de P.Vogel montrent que polyn{\^o}me $P$ ci-dessus,
vu comme {\'e}l{\'e}ment de $S_0$ s'envoie sur $0$ dans $\L_0$. La question de savoir si
le polyn{\^o}me $Q$ est nul dans $\L_0$ reste ouverte. Cette question est reli{\'e}e {\`a}
la conjecture de P. Deligne (cf \cite{D}, \cite{CM} et \cite{DM}) sur
l'existence d'une cat{\'e}gorie mono{\"\i}dale, lin{\'e}aire {\`a} coefficients polynomiaux,
universelle pour la famille des (super)alg{\`e}bres de Lie exceptionnelles.

\section{D{\'e}monstration de l'existence de $\X$}
\subsection{Les caract{\`e}res fondamentaux}
Dans \cite{Vo2}, huit caract{\`e}res gradu{\'e}s fondamentaux sont construits et tout
caract{\`e}re provenant d'une des superalg{\`e}bres de Lie simple mentionn{\'e}e se d{\'e}duit
de l'un d'eux; ces huit caract{\`e}res sont:
\begin{itemize}
\item $\X_{\sll}\co \L\go \QQ[\delta,\beta]$ (degr{\'e}($\delta$)$=1$,
  degr{\'e}($\beta$)$=2$).
\item $\X_{\osp}\co \L\go \QQ[\delta,\alpha]$ (degr{\'e}($\delta$)$=1$,
  degr{\'e}($\alpha$)$=1$).
\item $\X_{\DD}\co \L\go \QQ[\sigma_2,\sigma_3]$ (degr{\'e}($\sigma_2$)$=2$,
  degr{\'e}($\sigma_3$)$=3$).
\item les cinq caract{\`e}res $\X_L$ {\`a} valeurs dans $\QQ[t]$ pour les alg{\`e}bres de
  Lie exceptionnelles $L\in\{\gg_2,\,{\ff_4},\,\ee_6,\,\ee_7,\,\ee_8\}$.
\end{itemize}
On construit des caract{\`e}res {\`a} valeurs dans des quotients de $S$ {\`a} l'aide de ces
huit caract{\`e}res.

\subsection{Le caract{\`e}re $\X_1$ pour les familles $\sll$ et $\osp$}
Dans cette section, on construit un caract{\`e}re $\X_1$ qui r{\'e}alise le th{\'e}or{\`e}me
pour tous les triplets annul{\'e}s par les polyn{\^o}mes $P_\sll$ et $P_\osp$.
\subsubsection{Le foncteur $\Phi_\M$; $\X_{\sll(\infty)}=\X_{\osp(\infty)}$}
Les caract{\`e}res $\X_\sll$ et $\X_\osp$ permettent de d{\'e}finir des caract{\`e}res
$\X_{\sll(E)}$ et $\X_{\osp(E)}$ sur $\L$ pour une valeur formelle de la
superdimension de $E$. Dans les deux cas, en faisant tendre cette superdimension
vers $+\infty$, on obtient deux caract{\`e}res limites correspondant au m{\^e}me triplet
$(1,0,0)$. Cette premi{\`e}re construction permet de montrer la co{\"\i}ncidence de ces
deux caract{\`e}res limites.  

En reprenant les constructions d'{\'e}paicissement de \cite{BN}, il est facile de
construire une cat{\'e}gorie $\M$ mono{\"\i}dale, $\QQ$--lin{\'e}aire et un foncteur\\
$\Phi_\M\co \D\go\M$ qui factorise tous les foncteurs $\Phi_{\gl(E)}$ et
$\Phi_{\osp(E)}$ (La construction repose sur les repr{\'e}sentations standards de
$\gl$ et $\osp$): 

\rk{D{\'e}finitions et notations}

On d{\'e}signe les {\'e}l{\'e}ments de $\S_n$ comme produit de cycles disjoints. Par
exemple, $(1,2)(3)$ d{\'e}signe la transposition de $\S_3$ qui {\'e}change $1$ et
$2$. Si $\sigma\in\S_n$ est l'{\'e}l{\'e}ment
$(i^1_{1},\ldots,i^1_{k_1})\ldots(i^p_{1},\ldots,i^p_{k_p})$, on d{\'e}signe par
$<\sigma>$ le diagramme suivant de $\A((S^1)^{\amalg p},[n])$:
$$\input{cr_sigma.tex}$$
On note $\Delta$ le $(S^1,\emptyset)$--diagramme form{\'e} du seul cercle et on
d{\'e}finit $\Sigma_n$ (respectivement $\Sigma'_n$) comme le sous
$\QQ[\Delta]$--module libre de $\D_\Gamma([0],[n])$ de base
$\{<\sigma>,\,\sigma\in\S_n\}$ (respectivement
$\{<\sigma>,\,\sigma\in\S_n,\,\sigma\hbox{ a un point fixe}\}$).\\
Ensuite, on d{\'e}finit les cat{\'e}gories quotientes $\D_\osp$, $\D_\gl$ et $\D_\sll$ de\label{DM}
$\D_\Gamma$ (et les foncteurs quotients $\Phi_\osp$, $\Phi_\gl$ et $\Phi_\sll$)
obtenues en annulant les morphismes suivants:
\begin{itemize}
\item Pour $\D_\osp$:
        $$\input{rel_osp1.tex}$$
(On oublie l'orientation des courbes des diagrammes.)
        $$\input{rel_osp2.tex}$$
\item Pour $\D_\gl$:
        $$\begin{picture}(0,0)%
\epsfig{file=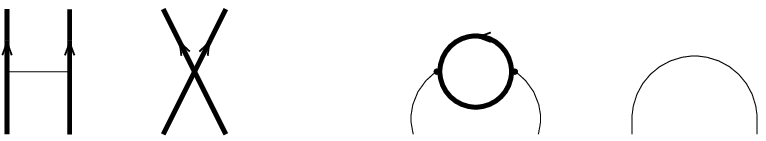}%
\end{picture}%
\setlength{\unitlength}{1973sp}%
\begingroup\makeatletter\ifx\SetFigFont\undefined%
\gdef\SetFigFont#1#2#3#4#5{%
  \reset@font\fontsize{#1}{#2pt}%
  \fontfamily{#3}\fontseries{#4}\fontshape{#5}%
  \selectfont}%
\fi\endgroup%
\begin{picture}(7256,1289)(558,-1007)
\put(6151,-661){\makebox(0,0)[b]{\smash{\SetFigFont{12}{14.4}{\rmdefault}{\mddefault}{\updefault}$\equiv$}}}
\put(1651,-511){\makebox(0,0)[b]{\smash{\SetFigFont{12}{14.4}{\rmdefault}{\mddefault}{\updefault}$\equiv$}}}
\end{picture}
$$
\item Pour $\D_\sll$:
        $$\hspace{-0.9cm}{\input{rel_sl.tex}}$$
Cette derni{\`e}re cat{\'e}gorie est en fait un quotient de
$\D_\Gamma\otimes\QQ[\Delta,\Delta^{-1}]$
\end{itemize}
\begin{prop}\label{gl,sl et osp}
Soit $E$ un superespace vectoriel de dimension finie. On adopte le choix
suivant pour la forme bilin{\'e}aire de $\osp(E)$ (respectivement $\gl(E)$ et
$\sll(E)$): $<x,y>=\frac12\hbox{str}_E(xy)$ (respectivement
$<x,y>=\hbox{str}_E(xy)$).
\begin{itemize}
\item Le foncteur $\Phi_{\osp(E),E}$ se factorise par $\Phi_\osp$.
\item Le foncteur $\Phi_{\gl(E),E}$ se factorise par $\Phi_\gl$.
\item Si la superdimension de $E$ est non nulle, le foncteur $\Phi_{\sll(E),E}$ se factorise par
        $\Phi_\sll$.  
\item On a les isomorphismes naturels:\\
        $\D_\gl([p],[q])\simeq \Sigma_{p+q}$ et $\D_\sll([p],[q])\simeq
        ({\Sigma_{p+q}}_{/\Sigma'_{p+q}})\otimes\QQ[\Delta,\Delta^{-1}]$.
\item $\forall K\in F_n$, $\exists! x\in\Sigma_n$ tel que
        $\Phi_\gl(K)=\Phi_\gl(x)$ et $\Phi_\sll(K)=\Phi_\sll(x)$.
\item Les caract{\`e}res $\X_\sll$ et $\X_\osp$ sont d{\'e}termin{\'e}s par:\\ 
        Si $K$ appartient {\`a} $F_n$, si $u$ est un {\'e}l{\'e}ment de $\L$, si les
        polyn{\^o}mes $P$ et $Q$ v{\'e}rifient $\X_\sll(u)=P(\delta,\beta)$ et
        $\X_\osp(u)=Q(\delta,\alpha)$ alors
        $$\Phi_\gl(u.K)=P(\Delta,1)\Phi_\gl(K)\quad\hbox{et}\quad
        \Phi_\osp(u.K)=Q(\Delta,1)\Phi_\osp(K)$$
\end{itemize}
\end{prop}
Pour la d{\'e}monstration de cette proposition, nous renvoyons aux arguments de
\cite{Vo2} sections 6.3 {\`a} 6.7.\\

Le concept de surface marqu{\'e}e et l'application d'{\'e}paicissement des diagrammes introduits par Bar-Natan (\cite{BN}), sont ici utilis{\'e}s pour construire $\M$: 

\rk{D{\'e}finition de $\M$}

On introduit d'abord la notion de $X$--surface marqu{\'e}e qui sera la donn{\'e}e d'une
classe d'isomorphisme de surface compacte {\`a} bord munie d'une bijection entre
l'ensemble fini $X$ et des tangentes non nulles au bord, prises en des points
distincts. Deux telles surfaces sont isomorphes s'il existe un diff{\'e}omorphisme
entre elles qui respecte les bijections de $X$ vers les tangentes de chacune.\\
On note $M(X)$ \label{MX}(respectivement $M[n]$) le $\QQ$--espace vectoriel engendr{\'e} par
les $X$--surfaces marqu{\'e}es (respectivement les $[n]$--surfaces marqu{\'e}es) quotient{\'e}
par les relations suivantes: Si $V$ est une $X$--surface marqu{\'e}e et si $V'$ est
obtenue en rempla\c{c}ant l'une des tangentes de $V$ par son oppos{\'e}e, alors:
$$V'\equiv-V$$ 
Comme pour $F(X)$, le groupe sym{\'e}trique $\S(X)$ agit sur $M[X]$ et en
reproduisant la construction de $\D$, on d{\'e}finit $\M$ par:
$$\hbox{Obj}(\M)=\{[n],\,n\in\NN\}$$
$$\M([p],[q])=M([p]\amalg[q])\simeq M[p+q]$$ 
Le produit tensoriel sur $\M$ est donn{\'e} par la r{\'e}union disjointe, la composition
$V_1\circ V_2$ est construite de la mani{\`e}re suivante: On joint chaque tangente
du but de $V_1$ {\`a} la tangente correspondante de la source de $V_2$ puis on
{\'e}paissit en appliquant la r{\`e}gle suivante:
$${\epsfbox{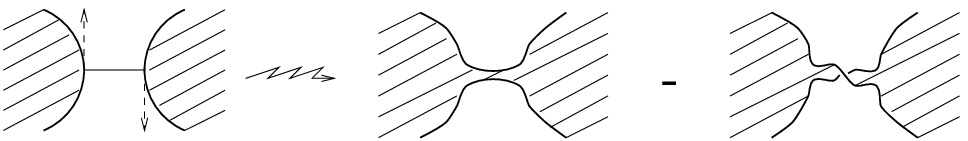}}$$  
Enfin $\Phi_\M$ prend les valeurs suivantes:
$$\begin{array}{cccccccc}\Phi_\M(&\put(-11,-8) {\epsfbox{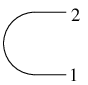}}&\,)=\frac12&
\put(-8,-8) {\epsfbox{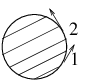}} & \qquad\Phi_\M( &
\put(-12,-8) {\epsfbox{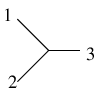}} & \,\,)= &\put(-8,-8) {\epsfbox{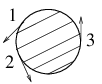}}
\end{array}$$ 
Il est facile de voir que dans $M[2]$ et $M[3]$, ces deux surfaces pr{\'e}sentent
les m{\^e}mes sym{\'e}tries que les diagrammes correspondants et un simple calcul montre
que l'image par $\Phi_\M$ d'une relation (IHX) est bien nulle.\\

On d{\'e}finit alors les foncteurs $\partial_{\gl}$ et $\partial_{\osp}$ entre les
cat{\'e}gories $\M$ et $\D_\Gamma$:\\ 
Si $V$ est une $[n]$--surface marqu{\'e}e ayant $p$ composantes de bord et si
$\epsilon$ est une orientation de $\partial V$, on d{\'e}signe par
$\partial(V,\epsilon)$ le diagramme form{\'e} des $p$ cercles orient{\'e}s constituant
le bord de $\partial V$ sur lequel on a fix{\'e} $n$ ``jambes'' num{\'e}rot{\'e}es aux lieux
des $n$ tangentes:
$$\input{bord.tex}$$
On pose alors:
$$\partial_{\gl}(V)=\sum_{x\hbox{ orientation de }V}\partial(V,\epsilon(x))$$
$$\partial_{\osp}(V)=2^{-p}\sum_{\epsilon\hbox{ orientation de }\partial V}
\partial(V,\epsilon)$$
(On a not{\'e} $\epsilon(x)$ l'orientation de $\partial V$ induite par $x$ et $p$ le
nombre de composante de bord de $V$.)
\begin{prop}\label{M=gl+osp} Soit $K\in\D([p],[q])$ et $V=\Phi_\M(K)$.\\
  $\partial_\gl(V)$ et $\partial_\osp(V)$ sont en fait des {\'e}l{\'e}ments de
  $\Sigma_{p+q}$ et on a:
  $$\Phi_\gl\circ\partial_{\gl}(V)=\Phi_\gl(K)\qquad
  \Phi_{\osp}\circ\partial_{\osp}(V)=\Phi_{\osp}(K)$$
\end{prop}
La d{\'e}monstration de cette propri{\'e}t{\'e} est la m{\^e}me que celle faite dans \cite{BN} pour
justifier la construction des applications d'{\'e}paicissement des diagrammes.\\

Le $\QQ$--espace vectoriel $M_c[0]$ (``c'' pour engendr{\'e} par les surfaces
connexes) est naturellement muni d'une structure d'alg{\`e}bre en prenant pour le
produit de deux surfaces leur somme connexe. Ainsi $M_c[0]$ est l'anneau\label{anns}\\
$\s\simeq\QQ[\delta,\alpha,\beta]_{/(\alpha\beta-\alpha^3)}$ o{\`u} $\alpha$,
$\beta$ et $\delta$ sont respectivement les classes de diff{\'e}o\-mor\-phis\-mes du
plan projectif r{\'e}el $\RR P^2$, du tore $S^1\times S^1$ et du disque $D^2$.\\
De plus, toujours par la somme connexe, les modules $M_c[n]$ sont munis d'une
structure de $\s$--module gradu{\'e} de type fini. Il est possible de d{\'e}crire une
famille g{\'e}n{\'e}ratrice de $M_c[n]$ de la mani{\`e}re suivante: Notons $\Sigma$ la
sph{\`e}re orient{\'e}e de $\RR^3$ et $D$ le disque ouvert de $\RR^2$. L'ensemble\\
$\E^p_n=\{V\simeq\Sigma\setminus D^{\amalg p}$ marqu{\'e}e par $[n]$ tangentes
respectant l'orientation, chaque composante connexe de $\partial V$ {\'e}tant munie
d'au moins une marque $\}$\\ est naturellement en bijection avec\\
$\{\sigma\in\S_n$ telles que $\sigma$ partitionne $[n]$ en $p$ orbite$\}$.\\
On notera $[\sigma]$ la classe dans $M_c[n]$ de la surface de $\E_n=\bigcup_p
\E^p_n$ correspondant {\`a} $\sigma\in\S_n$.\\ 
Il est facile de voir que $M_c[n]$ est engendr{\'e} par les surfaces marqu{\'e}es de
$\E_n$. De plus, {\`a} l'aide d'un diff{\'e}omorphisme de $\Sigma\setminus D^{\amalg p}$
renversant l'orientation, on montre l'indentit{\'e}: $[\sigma]=(-1)^n[\sigma^{-1}]$.\\

Si $L$ est l'alg{\`e}bre $\sll(E)$, on d{\'e}signe par $g_L$ le morphisme
d'alg{\`e}bre de $\s$ dans $\QQ[x]$ envoyant $(\alpha,\beta,\delta)$ sur
respectivement $(0,x^2,\sdim(E)x)$.\\ 
Si $L$ est l'alg{\`e}bre $\osp(E)$, on d{\'e}signe par $g_L$ le morphisme d'alg{\`e}bre de
$\s$ dans $\QQ[x]$ envoyant $(\alpha,\beta,\delta)$ sur respectivement
$(x,x^2,\sdim(E)x)$.\\

\begin{prop} Il existe un caract{\`e}re gradu{\'e} $\X_\M\co\L\go\s$ tel que, quels que
  soient $K\in F_n$ et $u\in\L$, on ait la propri{\'e}t{\'e}
  $$\Phi_\M(u.K)=\X_\M(u).\Phi_\M(K).$$ De plus $\X_\M$ mod $(\alpha)=\X_{\sll}$
  et $\X_\M$ mod $(\beta-\alpha^2)=\X_{\osp}$.\\
  En particulier pour $L=\sll(E)$ ou $L=\osp(E)$, on a:
  $$\X_{L}=g_L\circ\X_\M$$
\end{prop}
\begin{dem}
Pour commencer, $\Phi_\M(\L)$ est inclus dans la partie totalement
antisym{\'e}trique de $M_c[3]$ qui est le module libre de rang un engendr{\'e} par le
disque $[(1,2,3)]$ ($(1,2,3)$ d{\'e}signe le $3$--cycle de $\S_3$ qui envoie $1$ sur
$2$). Si $u$ est un {\'e}l{\'e}ment de $\L$, on peut donc d{\'e}finir $\X_\M(u)$ comme
l'{\'e}l{\'e}ment de $\s$ v{\'e}rifiant $\Phi_\M(u)=\X_\M(u).[(1,2,3)]$. Comme tout
diagramme $K$ de $F_n$ vu comme {\'e}l{\'e}ment de $\D([0],[n])$ peut s'{\'e}crire comme la
compos{\'e}e de l'{\'e}l{\'e}ment unit{\'e} de $\L$ ($u_0$ vu comme {\'e}l{\'e}ment de $\D([0],[3])$) et
d'un diagramme $\hat{K}$ de $\D([3],[n])$, on a:\\ $\Phi_\M(u.K)=
\Phi_\M(u\circ\hat{K})= \X_\M(u)\Phi_\M(u_0)\circ\Phi_\M(\hat{K})=
\X_\M(u)\Phi_\M(u_0\circ\hat{K})= \X_\M(u)\Phi_\M(K)$. Les relations entre
$\X_\M$, $\X_{\sll}$ et $\X_{\osp}$ sont des cons{\'e}quences directes de la
proposition \ref{M=gl+osp}.
\end{dem} 

\subsubsection{Construction du caract{\`e}re $\X_1$ pour les familles $\sll$ et
$\osp$} \label{Xslosp}
Le premier pas vers la construction de $\X$ consiste {\`a} modifier le caract{\`e}re
$\X_\M$ en un caract{\`e}re $\X_1$ {\`a} valeurs dans un quotient de $S$.
\begin{prop}
  Il existe des caract{\`e}res gradu{\'e}s
$$\X_1\co \L\go S_{/(P_{\sll}P_{\osp})}\qquad g\co S_{/(P_{\sll}P_{\osp})}\go\s$$ 
tels que $g\circ\X_1=\X_\M$ et pour toute superalg{\`e}bre de Lie $L$ de
type $\sll$ ou $\osp$ on ait: $$\X_L=g_L\circ\X_\M=f_L\circ\X_1$$
\end{prop}
Le reste de cette section est consacr{\'e} {\`a} la d{\'e}monstration de cette proposition.\\ 
Pour construire $g$, on peut calculer l'image par $\Phi_\M$ de la relation
(\ref{psi3}) ou utiliser la formule {\'e}tablie dans \cite{Vo2} pour chacun des
caract{\`e}res $\X_L$ cit{\'e}s:
\begin{equation}\label{somme des xn}
\sum_{n=0}^{\infty}\X_L(x_n)=
\frac{4tv+2t^2u-2t^4-3v-tu+7t^3-7t^2+2t}{(1-t-2u-2v) (1-t) (1-2t)}
\end{equation}
Dans les deux cas, on est amen{\'e} {\`a} d{\'e}finir $f\co S\go\s$ par:
$$f(t)=\delta-2\alpha$$
$$f(u)=2\beta+6\alpha^2-\alpha\delta$$
$$f(v)=16\alpha^3-2\delta\beta-2\delta\alpha^2$$
On trouve alors que le noyau de $f$ est l'id{\'e}al principal $(P_{\sll}P_{\osp})$
ce qui permet de factoriser $f$ en un morphisme d'alg{\`e}bre
$g\co S_{/(P_{\sll}P_{\osp})}\inj\s$. Pour pouvoir factoriser $\X_\M$ par $g$, il
faut montrer que $\X_\M(\L)$ est inclus dans Im$(g)$.\\
Pour cela on consid{\`e}re les superalg{\`e}bres de Lie suivantes qui admettent
plusieurs repr{\'e}sentations standards:
\begin{itemize} 
\item $L_1={\mathfrak{so}}(5)\simeq L'_1={\mathfrak{sp}}(4)$
\item $L_2=\sll(4)\simeq L'_2={\mathfrak{so}}(6)$
\item $L_3=\sll(2,1)\simeq L'_3=\osp(2,2)$
\end{itemize}
En suivant pour le Casimir de chacune de ces alg{\`e}bres les conventions de la
section pr{\'e}c{\'e}dente, chacun des isomorphismes envoie le Casimir $\Omega_i$ de
$L_i$ sur $\lambda_i\Omega'_i$ o{\`u} $\Omega'_i$ d{\'e}signe le Casimir de $L'_i$. En
particulier, ces alg{\`e}bres {\'e}tant isomorphes, on a: $\forall u\in\L,\,
\X_{L_i}(u)_{|x=1}=\X_{L'_i}(u)_{|x=\lambda_i}$. Ces paires de caract{\`e}res
pouvant {\^e}tre calcul{\'e}es par l'interm{\'e}diaire de $\X_\M$, cela fournit les
informations suivantes:\\
On d{\'e}finit six caract{\`e}res sur $\s$ par leurs valeurs sur le triplet
$(\delta,\alpha,\beta)$:
\begin{itemize}
\item $\Psi_1\co (\delta,\alpha,\beta) \im(10,2,4)$
\item $\Psi'_1\co (\delta,\alpha,\beta)\im(4,-1,1)$
\item $\Psi_2\co (\delta,\alpha,\beta) \im(4,0,1)$
\item $\Psi'_2\co (\delta,\alpha,\beta)\im(6,1,1)$
\item $\Psi_3\co (\delta,\alpha,\beta) \im(2,0,4)$
\item $\Psi'_3\co (\delta,\alpha,\beta)\im(0,-1,1)$
\end{itemize}
La compos{\'e}e de $\X_\M$ avec chacun de ces morphismes donne une renormalisation du
caract{\`e}re de l'une des six alg{\`e}bres ci-dessus et on a
$$\X_\M(\L)\subset\bigcap_{i=1}^{3}\Ker(\Psi'_i-\Psi_i)$$
De plus, ces informations sur l'image de $\X_L$ sont suffisantes pour d{\'e}montrer
la proposition car:
\begin{lem}
  $$\Im g=\bigcap_{i=1}^{3}\Ker(\Psi'_i-\Psi_i)$$
\end{lem}
\begin{dem}
  Pour voir que $\Im g\subset\bigcap_{i=1}^{3}\Ker(\Psi'_i-\Psi_i)$, il suffit de
  v{\'e}rifier que $\Psi'_i(x)$ et $\Psi_i(x)$ co{\"\i}ncident pour $i\in\{1,2,3\}$ et
  $x\in\{t,u,v\}$.  L'application $g$ {\'e}tant gradu{\'e}e et injective, pour montrer
  l'{\'e}galit{\'e} entre les espaces ci-dessus, il suffit de montrer l'{\'e}galit{\'e} de leurs
  dimensions en chaque degr{\'e}.  On fait le calcul pour $\Im g$ en utilisant les
  s{\'e}ries g{\'e}n{\'e}ratrices:
  $$\sum_{n=0}^\infty \dim((S)_n) t^n=\frac1{(1-t) (1-t^2) (1-t^3)}$$
  $$\sum_{n=0}^\infty \dim((S_{/(P_{\sll}P_{\osp})})_n) t^n=\frac{1-t^9} {(1-t)
    (1-t^2) (1-t^3)}$$
  $$\sum_{n=0}^\infty \dim((\s)_n) t^n=\frac{1-t^3}{(1-t) (1-t) (1-t^2)}$$
  $$\sum_{n=0}^\infty \hbox{codim}((\Im g)_n) t^n =$$
  $$\frac{1-t^9}{(1-t)(1-t^2)(1-t^3)}-\frac{1-t^3}{(1-t)(1-t)(1-t^2)}=
  \frac3{1-t}-3-2t-t^2-t^3$$
  Donc l'image de g est de codimension $3$ en chaque degr{\'e} $n\geq 4$.\\
  D'autre part un calcul facile donne en degr{\'e} $n\geq 3$:
  $$\hbox{det}(\left((\Psi'_i-\Psi_i)
    (\delta^{n+1-j}\alpha^{j-1})\right)_{i,j=1\ldots
    3})=2^{n+1}6^{n-2}(5\times4^{n-2}-2\times10^{n-2})$$
  Ce d{\'e}terminant ne s'annule que pour $n=3$, ce qui assure qu'en chaque degr{\'e}
  $n\geq 4$, les formes lin{\'e}aires induites par $(\Psi'_i-\Psi_i)$ sont
  ind{\'e}pendantes et donc que l'intersection de leur noyau est de codimension
  trois. Il reste {\`a} v{\'e}rifier le lemme en bas degr{\'e}s, ce qui peut se faire
  directement, terminant ainsi la d{\'e}monstration.
\end{dem} 

\paragraph{\bf Remarque}
Il semble difficile de traduire de mani{\`e}re g{\'e}om{\'e}trique la restriction de l'image
vde $\X_\M$ bien que ce caract{\`e}re ait une construction par {\'e}pai\-cis\-se\-ment des
diagrammes.

\subsection{Le caract{\`e}re $\X_2$ pour les familles $\sll$, $\osp$ et
$\DD$}\label{car_X2} En notant $R=\QQ[a,b,c]_{/(a+b+c)}$ et $k$ son corps de
fractions, il existe une $R$--superalg{\`e}bre de Lie $\DD$ munie d'un {\'e}l{\'e}ment de
Casimir $\Omega\in \DD\otimes \DD$ telle que \label{D21}
\begin{itemize}
\item En {\'e}tendant les coefficients {\`a} $k$, $\Omega$ est un {\'e}l{\'e}ment de Casimir non
  d{\'e}g{\'e}n{\'e}r{\'e} de $\wb\DD\otimes_k \wb\DD$ (o{\`u}
  $\wb\DD=\DD\otimes_Rk$).
\item Vue comme sous-alg{\`e}bre de $\wb\DD$, $\DD$ se d{\'e}compose en somme directe
  orthogonale: $\DD=L_1\oplus L_2\oplus L_3\oplus X$ o{\`u} $L_i$ est une sous
  alg{\`e}bre isomorphe {\`a} $\sll_2(\QQ)\otimes R$, $X$ est la partie impaire de $\DD$
  et l'action de $L_1\oplus L_2\oplus L_3$ sur $X$ induite par le crochet est
  donn{\'e}e par le produit tensoriel des repr{\'e}sentations standards de $L_1$, $L_2$
  et $L_3$.\\La d{\'e}composition {\'e}tant orthogonale, on a:
  $\Omega=-a\omega_1-b\omega_2-c\omega_3+\pi$ avec $\omega_i\in L_i\otimes L_i$
  et $\pi\in X\otimes X$.
\item Le groupe sym{\'e}trique $\S_3$ agit sur $R$ (par permutation de $\{a,b,c\}$)
  et sur $\DD$ (par permutation des trois copies de $\sll_2$) de mani{\`e}re
  compatible avec la structure de $R$--module. Cette action laisse le
  Casimir invariant et commute au crochet de Lie.
\item Le foncteur $\Phi_{\DD}\co \D_b\go\Mod_{\DD}$ induit un caract{\`e}re gradu{\'e}\\ 
  $\Xd\co \L\go\QQ[\sigma_2,\sigma_3]$ (o{\`u} $\QQ[\sigma_2,\sigma_3]$ est le sous
  anneau de $R$ engendr{\'e} par $\sigma_2=ab+bc+ca$ de degr{\'e} $2$ et $\sigma_3=abc$
  de degr{\'e} $3$).
\item Soit $\alpha$ un nombre complexe et $\phi\co R\go\CC$ est le morphisme
d'anneau d{\'e}fini par:\\ $\phi(a)=\alpha$, $\phi(b)=1$ et $\phi(c)=-1-\alpha$.
Alors $\DD\otimes_\phi\CC\simeq {\mathfrak D}_{2\,1,\alpha}$ (not{\'e} parfois
$\osp_\alpha(4,2)$) et $\X_{{\mathfrak D}_{2\,1,\alpha}}=\phi\circ\Xd$.
\end{itemize}
Tous ces r{\'e}sultats sont imm{\'e}diatement cons{\'e}quences de ceux {\'e}tablis dans \cite{Vo2}
sections 6.10 {\`a} 6.15 (o{\`u} l'anneau de coefficients consid{\'e}r{\'e} est
$\ZZ[a,b,c]_{/(a+b+c)}$ au lieu de $R$). De plus, il y est d{\'e}montr{\'e} que
l'op{\'e}rateur $\Psi$ a pour valeurs propres $\{2a,2b,2c\}$. Ainsi le changement de
variables $\sigma_2\mapsto-\frac{u}2$ et $\sigma_3\mapsto\frac{v}4$ permet de
consid{\'e}rer $\Xd$ comme un caract{\`e}re {\`a} valeurs dans $S_{/(t)}$.\\
\begin{prop} Il existe un caract{\`e}re gradu{\'e} $\X_2\co \L\go S_{/(tP_{\sll}P_{\osp})}$
  factorisant les caract{\`e}res $\X_1$ et $\Xd$\\ (c'est {\`a} dire: $\X_2$ mod
  $(P_{\sll}P_{\osp})=\X_1$ et $\X_2$ mod $(t)=\Xd$).
\end{prop}
\dem
  L'anneau $S$ {\'e}tant factoriel, l'intersection des id{\'e}aux $(t)$ et
  $(P_{\sll}P_{\osp})$ est l'id{\'e}al $(tP_{\sll}P_{\osp})$. Donc le diagramme
  commutatif suivant est un carr{\'e} cart{\'e}sien:
\begin{eqnarray*}
  S_{/(tP_{\sll}P_{\osp})}&\go&S_{/(P_{\sll}P_{\osp})}\\ 
  \downarrow&&\downarrow\\ 
  S_{/(t)}&\go&S_{/(t)+(P_{\sll}P_{\osp})}
\end{eqnarray*}
  Ainsi, il suffit de montrer que le diagramme suivant est commutatif:
\begin{eqnarray*}
  \L&\stackrel{\X_1}{\go}&S_{/(P_{\sll}P_{\osp})}\\ 
  \Xd\downarrow&&\downarrow\\ 
  S_{/(t)}&\go&S_{/(t)+(P_{\sll}P_{\osp})}\inj S_{/(t)+(v)}\times
  S_{/(t)+(27v^2-8u^3)}
\end{eqnarray*}
  (On a rajout{\'e} la derni{\`e}re injection pour aboutir dans un produit d'anneaux
  int{\`e}gres).
\begin{itemize}
\item La co{\"\i}ncidence des deux caract{\`e}res $\L\go S_{/(t)+(27v^2-8u^3)}$ est une
cons{\'e}\-quence directe de l'isomorphisme entre les superalg{\`e}bres de Lie\\
$L={\mathfrak D}_{2\,1,\alpha}\simeq \osp(4,2)$ pour
$\alpha\in\{1,-2,-\frac12\}$ (on a alors par exemple $(a,b,c)=(1,-2,1)$ donc
$t(L)=0$, $u(L)=6$ et $v(L)=-8$). Le premier caract{\`e}re correspond {\`a} factoriser
$\X_L$ par $\Xd$, l'autre revient {\`a} factoriser $\X_L$ par $\X_1$ ou $\X_\M$.
\item Le quotient de $\Xd$ par l'id{\'e}al $(v)$ correspond au caract{\`e}re induit par
  l'alg{\`e}bre de Lie $L={\mathfrak D}_{2\,1,\alpha}$ dans le cas d{\'e}g{\'e}n{\'e}r{\'e} o{\`u}
  $\alpha$ vaut $0$ ou $-1$. Dans ce cas, $L$ contient un id{\'e}al $\h\simeq
  \psl(2,2)$ et le Casimir $\Omega$ appartient {\`a} $\h\otimes\h$. Ainsi une
  application directe du lemme \ref{Phi_b} pour $(\phi,\h\inj L)$, montre que
  $\X_L$ co{\"\i}ncide avec $\X_{\psl(2,2)}$ qui se factorise par $\X_1$. Ce
  caract{\`e}re est en fait le caract{\`e}re augmentation (nul sauf en degr{\'e} 0) de $\L$.\endproof
\end{itemize}

\subsection{Le caract{\`e}re $\X_3$ pour les familles $\sll$, $\osp$, $\DD$ et
  $\sll_2$}
\subsubsection{Th{\'e}or{\`e}me d'existence}
Le cas de l'alg{\`e}bre $\sll_2$ est un cas particulier: On peut choisir pour
  $(t,u,v)$ n'importe quelle valeur pourvu que le polyn{\^o}me
  $P_{\sll_2}=v-ut+t^3$ soit nul; de plus ce choix n'affecte bien s{\^u}r pas la
  valeur de $\X_{\sll_2}\in\QQ[x]$. On peut donc consid{\'e}rer $\X_{\sll_2}$ comme un
  caract{\`e}re non surjectif {\`a} valeurs dans $S_{/P_{\sll_2}}$ en le composant avec
  l'inclusion de $\QQ[x]$ dans $S_{/P_{\sll_2}}$ qui envoie $\X_{\sll_2}(t)$ sur
  $t$.
\begin{prop} Il existe un caract{\`e}re gradu{\'e} $\X_3\co \L\go
  S_{/(tP_{\sll}P_{\osp}P_{\sll_2})}$ factorisant les caract{\`e}res $\X_2$ et
  $\X_{\sll_2}$\\
  (c'est {\`a} dire: $\X_3$ mod $(P_{\sll}P_{\osp}t)=\X_2$ et $\X_3$ mod
  $(P_{\sll_2})=\X_{\sll_2}$).
\end{prop}
Bien que le cas $\X_{\sll_2}$ ait d{\'e}j{\`a} {\'e}t{\'e} trait{\'e} pr{\'e}c{\'e}demment (par exemple
comme sous-cas de $\sll_n$), son {\'e}tude dans ce cadre permet d'affiner le
caract{\`e}re $\X_2$.
\begin{dem...}
De m{\^e}me que pr{\'e}c{\'e}demment, $S$ {\'e}tant factoriel, l'intersection des id{\'e}aux
$(P_{\sll_2})$ et $(P_{\sll}P_{\osp}t)$ est l'id{\'e}al
$(tP_{\sll}P_{\osp}P_{\sll_2})$. Donc le diagramme commutatif suivant est un
carr{\'e} cart{\'e}sien:
\begin{eqnarray*}
S_{/(tP_{\sll}P_{\osp}P_{\sll_2})}&\go&S_{/(tP_{\sll}P_{\osp})}\\
\downarrow&&\downarrow\\
S_{/(P_{\sll_2})}&\go&S_{/(P_{\sll_2})+(tP_{\sll}P_{\osp})}
\end{eqnarray*}
Donc il suffit de montrer que le diagramme suivant est commutatif:
\begin{eqnarray*}
\L&\stackrel{\X_2}{\go}&S_{/(tP_{\sll}P_{\osp})}\\
\X_{\sll_2}\downarrow&&\downarrow\\
\QQ[t]\subset S_{/(P_{\sll_2})}&\go&S_{((P_{\sll_2})+(tP_{\sll}P_{\osp})} \inj
S_{/I_1} \times S_{/I_2} \times S_{/I_3} \times S_{/I_4} \times S_{/I_5}\\ \\
\hbox{avec}&&I_1= (P_{\sll_2})+(2u-t^2) = (2v+t^3)+(2u-t^2)\\
&&I_2= (P_{\sll_2})+(u-5t^2)=(v-4t^3)+(u-5t^2)\\
&&I_3= (P_{\sll_2})+(8u-5t^2)=(8v+3t^3)+(8u-5t^2)\\
&&I_4= (P_{\sll_2})+(u-t^2)=(v)+(u-t^2)\\
&&I_5= (P_{\sll_2})+(t^2)=(v-ut)+(t^2)
\end{eqnarray*}
L'anneau $S_{((P_{\sll_2})+(tP_{\sll}P_{\osp})}$ s'injecte dans le produit
d'anneaux ci-dessus dont tous les facteurs, sauf le dernier, sont int{\`e}gres et en
fait isomorphes {\`a} $\QQ[t]$ ce qui permet de consid{\'e}rer la fl{\`e}che horizontale du
bas sur chacun de ces facteurs:\\
$\QQ[t]\subset
S_{/(P_{\sll_2})}\go S_{((P_{\sll_2})+(tP_{\sll}P_{\osp})} \go
S_{/I_k}\simeq\QQ[t])$\\
comme {\'e}tant l'identit{\'e} sur $\QQ[t]$.
\begin{itemize}
\item On a $I_1=I_{\sll_n}$ pour $n=2$ et la commutativit{\'e} du diagramme
  correspond sur ce facteur au fait que $\X_{\sll}$ factorise $\X_{\sll_2}$.
\item De m{\^e}me $I_2=I_{\so_3}$, $I_3=I_{\sp_2}$ et $I_4=I_{\so_4}$; Les
  isomorphismes bien connus entre les alg{\`e}bres de Lie
  $\so_3(\CC)\simeq\sp_2(\CC)\simeq\sll_2(\CC)$ et\\ 
  $\so_4(\CC)\simeq\sll_2(\CC)\oplus\sll_2(\CC)$ permettent facilement de voir
  que le diagramme commute sur chacun des facteurs correspondant.
\end{itemize}

La seule difficult{\'e} est de montrer la commutativit{\'e} du diagramme sur le dernier
facteur: l'image du caract{\`e}re $\X_{\sll_2}$ est en fait contenue dans $\QQ[t]$
et sa compos{\'e}e avec $S_{/(P_{\sll_2})}\go S_{/I_5}$ est donc nulle en degr{\'e}
sup{\'e}rieur ou {\'e}gal {\`a} deux (car $I_5$ contient l'id{\'e}al $(t^2)$). Il s'agit donc de
montrer que $\X_2$ mod $I_5$ est aussi nul en degr{\'e} sup{\'e}rieur o{\`u} {\'e}gal {\`a} deux.\\
On a d{\'e}j{\`a} vu que $\X_2$ mod $I_{\psl(2,2)}(=(v)+(t))$ {\'e}tait nul en degr{\'e}
strictement positif.  Mais $I_5\subset((v)+(t))$ et l'application quotient
$S_{/I_5}\go S_{/(v)+(t)}$ est un isomorphisme en degr{\'e} pair donc $\X_2$ mod
$I_5$ est nul en degr{\'e} pair.\\  
Si $x\in \L$ est de degr{\'e} $n=2p+1$, on {\'e}crit
$\X_2(x)=\lambda(x)u^pt+\lambda'(x)u^{p-1}v+r(x)$ avec
$\lambda(x),\lambda'(x)\in\QQ$ et $r(x)\in((v-tu)+(t^2))$.\\
On a alors: $\X_2(x)\equiv (\lambda(x)+\lambda'(x))u^pt\hbox{ mod
  }((v-tu)+(t^2))$. Il s'agit de montrer que $\lambda+\lambda'=0$. Revenant aux
  caract{\`e}res fondamentaux $\X_{\DD}$ {\`a} valeurs dans $\QQ[\sigma_2,\sigma_3]$ et
  $\X_\sll$ {\`a} valeurs dans $\QQ[\delta,\beta]$ cela revient {\`a} montrer que:\\
si $\X_{\DD}(x)\equiv \mu\sigma_3\sigma_2^{p-1}$ mod $(\sigma_3^2)$\\
alors $\X_\sll(x)\equiv -\mu\delta\beta^p$ mod $(\delta^2)$.\\
Ce point est plus d{\'e}licat et sa d{\'e}monstration fait l'objet des sections
suivantes.
\end{dem...}

\subsubsection{Propri{\'e}t{\'e}s du caract{\`e}re $\X_\DD$}\label{prop D21}
Le but de cette section est de d{\'e}montrer que la r{\'e}duction modulo $a^2$ de
$\X_\DD$ peut {\^e}tre calcul{\'e}e par l'interm{\'e}diaire du foncteur $\Phi_{\gl_{2,2}}$.
On reprend ici les notations introduites dans la partie \ref{car_X2}. De plus, on
note $F_n^+$ le sous-module de $F_n$ engendr{\'e} par les diagrammes ayant au moins
une boucle. Afin d'harmoniser les notations, on utilisera souvent $a_1$ pour $a$,
$a_2$ pour $b$ et $a_3$ pour $c$. On a alors

\begin{lem} \label{PhiDD}
Le foncteur $\Phi_{\DD}$ a les propri{\'e}t{\'e}s suivantes:\\
\begin{itemize}
\item $\Phi_{\DD}\co F_2\go\DD^{\otimes 2}$ est nulle
\item $\forall n\geq 3,\, \Phi_{\DD}\co F_n^+\go\DD^{\otimes n}$ est en
  fait {\`a} valeurs dans \\
  $\bigoplus E_1\otimes E_2\otimes\ldots\otimes E_n$ o{\`u} $E_i$ parcourt
  $\{a_1L_1,a_2L_2,a_3L_3,X\}$. En particulier, par restriction, la dualisation
  dans $\wb\DD$ induit: 
  $$F_n^+\go\Mod_\DD(\DD^{\otimes n},R)$$
  (morphismes de $\DD$--modules).
\end{itemize}
\end{lem}
\begin{dem}
Pour la premi{\`e}re affirmation, il suffit de consid{\'e}rer que $F_2$ est un
$\L$--module libre de rang $1$ engendr{\'e} par le diagramme 
$${\epsfbox{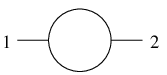}}$$
Dont l'image par $\Phi_{\DD}$ est nulle.\\

Pour d{\'e}montrer la seconde affirmation, on choisit un diagramme $K$ con\-nexe de
$F_n^+$ et on d{\'e}compose le calcul de $\Phi_{\DD}(K)$ de la mani{\`e}re suivante: on
colorie chaque ar{\^e}te de $K$ par $\Omega$ puis on utilise l'application de
r{\'e}duction $\varphi=\bigotimes_x \varphi_x$ o{\`u} $\varphi_x$ est ``la r{\'e}duction
autour du sommet trivalent $x$'' par l'application trilin{\'e}aire antisym{\'e}trique
invariante $<[\cdot,\cdot],\cdot>$ (ceci revient {\`a} d{\'e}composer $K$ comme compos{\'e}
du produit tensoriel d'autant d'{\'e}l{\'e}ments de Casimir que $K$ comporte d'ar{\^e}tes
avec un diagramme induit par une permutation puis avec un diagramme produit
tensoriel de $n$ fois l'identit{\'e} et d'autant de fois le diagramme $<[.,.],.>$
que $K$ comporte de sommets). On peut remarquer que l'application $\varphi$ est
bien d{\'e}finie d{\`e}s que les ar{\^e}tes de $K$ sont colori{\'e}es par une application de
l'ensemble de ses ar{\^e}tes vers $S^2L$ et on a aussit{\^o}t la formule consistant {\`a}
``d{\'e}velopper'' le diagramme $K$ colori{\'e} par $\Omega$:
$$\Phi_L(K)=\sum_{c\in{\cal C}} \varphi(c(K))$$
o{\`u} ${\cal C}=\{c:\{\hbox{ar{\^e}tes de
  }K\}\go\{-a_1\omega_1;-a_2\omega_2;-a_3\omega_3;\pi\}\}$ et $c(K)$ d{\'e}signe $K$
colori{\'e} par $c$.\\ 
Or pour l'alg{\`e}bre $\wb{\DD}$, l'application $<[\cdot,\cdot],\cdot>$ est nulle
sur $L_i\otimes L_j \otimes X$ si $i\neq j$ et sur $X\otimes X\otimes X$. Donc
$\varphi(c(K))$ est nul d{\`e}s que l'un des sommets de $c(K)$ n'est pas d'une des
deux formes suivantes:
$$\input{d21admis.tex}$$ 
Un diagramme colori{\'e} dont tous les sommets sont de type $(1)$ ou $(2)$ sera dit
admissible.\\ 
De plus, $<[\cdot,\cdot],\cdot>$ envoie $L_i\otimes L_i\otimes L_i$ sur
$a_i^{-1}R$ et $L_i\otimes X\otimes X$ sur $R$; il suffit donc de voir (par
exemple pour $i=1$) que si $c(K)$ est admissible, si $m$ d{\'e}signe le nombre
d'ar{\^e}tes int{\'e}rieures de $c(K)$ de couleur $-a_1\omega_1$, $s_3$ le nombre de
sommets de $c(K)$ de type $(1)$ pour $i=1$ et $s_1$ le nombre de sommets
univalents de $K$ dont l'ar{\^e}te issue (ar{\^e}te dite ext{\'e}rieure) est de couleur
$-a_1\omega_1$, alors $m\geq s_3+s_1$.\\ 
\begin{itemize}
\item Si $c(K)$ est un diagramme $K$ {\`a} $d$ boucles pour lequel toutes les ar{\^e}tes
  sont de couleur $-a_1\omega_1$, les formules reliant le nombre d'ar{\^e}tes et le
  nombre de sommets d'un diagramme {\`a} son degr{\'e} donnent: $m=2s_1+3d-3$ et
  $s_3=s_1+2d-2$ donc $m-s_3-s_1=d-1\geq0$ car $d>0$.
\item Sinon, le sous-graphe de $c(K)$ form{\'e} des ar{\^e}tes colori{\'e}es par
$-a_1\omega_1$ est form{\'e} de plusieurs graphes connexes dont les sommets
univalents proviennent soit des sommets univalents de $K$, soit de sommets de
type $(2)$ de $c(K)$. De plus, $K$ {\'e}tant connexe, chaque composante connexe $K_j$
de ce sous-graphe poss{\`e}de un nombre $s_{2,j}>0$ de sommets univalents provenant
de sommets de type $(2)$ de $c(K)$. On note de m{\^e}me $d_j$ le nombre de boucles
de $K_j$, $n_j$ son nombre de sommets univalents, $m_j$ son nombre d'ar{\^e}tes et
$s_{3,j}$ son nombre de sommets trivalents. On a $s_1=\sum_j (n_j-s_{2,j})$,
$s_3=\sum_j s_{3,j}$ et $m=\sum_j m_j$. Ceci donne l'in{\'e}galit{\'e} cherch{\'e}e en
sommant sur $j$ les formules $m_j=2n_j+3d_j-3$ et $s_{3,j}=n_j+2d_j-2$. On
obtient: $m-s_3-s_1=\sum_j(d_j+s_{2,j}-1)\geq0$ car $s_{2,j}>0$.
\end{itemize}
Ceci conclut la d{\'e}monstration du lemme.
\end{dem}

$$\begin{array}{|c|}\hline\\
\input{figd21.tex}\\\hline\end{array}$$

Consid{\'e}rons maintenant $K$ un diagramme repr{\'e}sentant un {\'e}l{\'e}ment de $F_6^+$ de
degr{\'e} impair. Soit $G$ le sous-groupe de $\S_6$ laissant fixe
$\{\{1;2\};\{3;4\};\{5;6\}\}$. On d{\'e}finit le morphisme de groupe $\epsilon$ sur
$G$ en remarquant que tout {\'e}l{\'e}ment $\sigma\in G$ agit par
$\epsilon(\sigma)\in\{-1;+1\}$ sur le diagramme $D_2$ ci-dessus de $F_6$.\\ 
On notera les {\'e}l{\'e}ments de $G$ comme produits de cycles disjoints; par exemple
on a $\epsilon((1,4,2,3))=1$. Soit 
$$\tilde{K}=\frac 1{48}\sum_{\sigma\in G}\epsilon(\sigma)\sigma(K).$$ Le
morphisme $D_1\in\D([6],[3])$ envoie $\tilde{K}$ sur $\K\in
F_3\otimes_{\S_3}\QQ^-\simeq\L$. De plus, le degr{\'e} de $\K$ dans $\L$ est {\'e}gal {\`a}
celui de $K$ moins deux.\\ 
Soit $\alpha \in\QQ$ tel que $\Xd(\K)\equiv\alpha\sigma_3\sigma_2^p$ mod
$(\sigma_3^2)\,\, \equiv\alpha a_1(a_2a_3)^{p+1}$ modulo $a_1^2$.\\ 

On construit un tenseur dont l'image par $\Phi_\DD(K)$ (vu comme {\'e}l{\'e}ment de
$\Hom(L^{\otimes 6},\QQ)$) est congrue {\`a} $\X_\DD(\tilde{K}_3)$ modulo
$(a_1^2)$:\\

Soit $x_0=-a_1\varphi(D_3)\in L_1^{\otimes 3}$ de sorte que consid{\'e}rant
$\Phi_\DD(\K)$ comme un {\'e}l{\'e}ment de $\Mod(L^{\otimes 3},\QQ)$, on ait: \\
$$\Phi_\DD(\tilde K_3)(x_0)=\X_\DD(\tilde K_3)\Phi_{\wb\DD}\left(\put(1,-5){\epsfbox{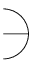}}\hspace{.3cm}\right)(x_0)=\frac{12}{a_1}
\X_\DD(\tilde K_3)$$
$$\equiv 12\alpha(a_2a_3)^{p+1}+a_1r\mbox{ avec }r\in R$$ 
D'autre part, le calcul montre que $\Phi_\DD(D_1') (x_0)$ s'{\'e}crit $x_1+a_1x_2$
avec $x_1\in X^{\otimes 6}$ et $x_2\in\DD^{\otimes 6}$, le tenseur $x_1$ ayant
de plus la propri{\'e}t{\'e} d'{\^e}tre $(G,\epsilon)$--invariant (i.e. $\forall\sigma\in
G,\,\sigma(x_1)=\epsilon(\sigma)x_1$). Ainsi
$$\Phi_\DD(\K) (x_0)=\Phi_\DD(\tilde{K})\circ\Phi_\DD(D_1') (x_0)=
\Phi_\DD(\tilde{K}) (x_1)+a_1\Phi_\DD(\tilde{K}) (x_2)$$
mais $\Phi_\DD(\tilde{K}) (x_2)$ est un {\'e}l{\'e}ment de $R$ d'apr{\`e}s le lemme
pr{\'e}c{\'e}dent et puisque $x_1$ est $(G,\epsilon)$--invariant, on a l'{\'e}galit{\'e}
$\Phi_\DD(\tilde{K}) (x_1)=\Phi_\DD(K) (x_1)$; donc
$$\Phi_\DD(K) (x_1)\equiv12\alpha(a_2a_3)^{p+1}\hbox{ modulo }(a_1).$$
On construit maintenant une forme lin{\'e}aire $\eta''$ sur $\gl_{2,2}^{\otimes 6}$
telle que l'application $K\im\alpha$ se factorise en
$\eta''\circ\Phi_{\gl_{2,2}}$.\\
Soit $I$ le noyau du morphisme d'anneau $f\co R\go\QQ$ d{\'e}fini par $f(a_1)=0$ et
$f(a_2)=1$. La r{\'e}duction des coefficients gr{\^a}ce {\`a} $f$ induit
$p\co \DD\go\DD'=\DD\otimes_f\QQ$. Le morphisme d'alg{\`e}bre $p$ est surjectif et
poss{\`e}de une section d'espace vectoriel $\ZZ_2$--gradu{\'e}: $s\co \DD'\inj\DD$ induite
par le monomorphisme $\QQ\inj R$ (on a $p\circ s=Id_{\DD'}$). On notera $X'$
l'image par $s$ de la partie impaire de $\DD'$ de sorte que $X$ se d{\'e}compose en
$X'\oplus (IX')$. Il se trouve que $<X',X'>$ est inclus dans $\QQ$ et $x_1$
appartient {\`a} ${X'}^{\otimes 6}$. Par cons{\'e}quent, les formes $\QQ$--lin{\'e}aires
$y\im f(<x_1,y>_{\wb\DD^{\otimes 6}})$ et $y\im <x_1,s\circ
p(y)>_{\wb\DD^{\otimes 6}}$ sont {\'e}gales. Mais consid{\'e}rant $\Phi_\DD(K)$ comme
un {\'e}l{\'e}ment de $\DD^{\otimes 6}$, on a:\\
$\begin{array}{lc}
&<x_1,\Phi_\DD(K)>_{\wb\DD}\equiv 12\alpha(a_2a_3)^{p+1}\hbox{ modulo }(a_1)\\
\hbox{donc }&f(<x_1,\Phi_\DD(K)>_{\wb\DD^{\otimes
6}})=(-1)^{p+1}12\alpha\\ \hbox{et donc }
&<x_1,(s\circ p)_*(\Phi_\DD(K))>_{\wb\DD^{\otimes 6}})=(-1)^{p+1}12\alpha.
\end{array}$\\
Ainsi l'application $\eta\co F_6\go\QQ$, qui associe au diagramme $K'$ l'{\'e}l{\'e}ment\\ 
$\frac{(-1)^{p+1}}{12} f(<x_1,\Phi_\DD(K')>_{\wb\DD^{\otimes 6}})$, se factorise
en $\eta=\eta'\circ\Phi_{\DD',p_*\Omega}$ o{\`u} $\eta'$ est une forme lin{\'e}aire sur
${\DD'}^{\otimes 6}$ nulle sur l'orthogonal de $((\DD')_{\bar 1})^{\otimes 6}$.\\
Consid{\'e}rons alors les morphismes de superalg{\`e}bres de Lie suivants:
\begin{eqnarray*}
  &\psl_{2,2}&\simeq\h\stackrel{i}{\inj}\DD'\\ 
  &\downarrow j\\ 
  \gl_{2,2}\stackrel{q}{\twoheadrightarrow} &\pgl_{2,2}&
\end{eqnarray*}
$\h$ est l'id{\'e}al de $\DD'$ d{\'e}j{\`a} rencontr{\'e} {\`a} la section \ref{car_X2}. Il est
facile de constater que si l'on muni $\gl_{2,2}$ du Casimir induit par la
repr{\'e}sentation standard, alors son image par $q_*$ dans $\pgl_{2,2}^{\otimes 2}$
appartient en fait {\`a} $j_*(\psl_{2,2}^{\otimes 2})$.\\
On peut maintenant appliquer le lemme \ref{PhiDD} {\`a} chacun de ces morphismes en
prenant pour $\phi$ l'identit{\'e} de $\QQ$:
$$\Phi_{\DD',p_*\Omega}(K)=i_*\circ\Phi_{\psl_{2,2}}(K)$$
$$q_*(\Phi_{\gl_{2,2}}(K))=\Phi_{\pgl_{2,2}}(K)=j_*(\Phi_{\psl_{2,2}}(K))$$ 
Or les morphismes $i$, $j$ et $q$ sont bijectifs en degr{\'e} impair. Comme $\eta'$ ne
d{\'e}pend que de la composante sur $((\psl_{2,2})_{\bar 1})^{\otimes 6}$ de
$\Phi_{\psl_{2,2}}(K)$, il existe une forme lin{\'e}aire $\eta''\co \gl_{2,2}^{\otimes
  6}\go\QQ$ nulle sur l'orthogonal de $((\gl_{2,2})_{\bar 1})^{\otimes 6}$
v{\'e}rifiant $\eta=\eta''\circ\Phi_{\gl_{2,2}}$.

\subsubsection{Propri{\'e}t{\'e}s du foncteur $\Phi_{\gl}$}
Dans toute cette section, $E$ d{\'e}signe un supermodule (c'est {\`a} dire un
$\QQ$--espace vectoriel muni d'une $\ZZ_2$--graduation) et on continue
d'identifier $\D_\gl([p],[q])$ avec $\Sigma_{p+q}$.\\
On notera $\D_{\gl_0}$ \label{Dgl0}la cat{\'e}gorie quotient{\'e}e de $\D_\gl$ par $\Delta\equiv0$
et $\Phi_{\gl_0}\co \D_\Gamma\go\D_{\gl_0}$.\\
On peut maintenant faire quelques remarques sur l'image de $\Phi_\gl$:
\begin{itemize}
\item $\partial_\gl$ est d{\'e}termin{\'e} sur $M[n]$ par la formule
  $$\forall\sigma\in\S_n,\,\partial_\gl([\sigma])=<\sigma>+(-1)^n<\sigma^{-1}>.$$
  En cons{\'e}quence, $\Phi_{\gl}(F_n)=\Phi_{\gl}\circ\partial_\gl\circ\Phi_\M(F_n)$
  est invariant par l'endomor\-phisme de $\Sigma_n$ d{\'e}fini sur sa base par\\
  $<\sigma>\im\frac12(<\sigma>+(-1)^n<\sigma^{-1}>)$.
\item Ensuite $\Phi_\M$ respecte le degr{\'e} et $\beta$ est de degr{\'e} pair donc
l'image par $\Phi_\M$ d'un {\'e}l{\'e}ment de degr{\'e} impair de $F_n$ est une combinaison
lin{\'e}aire {\`a} coefficients dans $\QQ[\beta]$ d'{\'e}l{\'e}ments
$[\sigma]\in\bigcup_p\E^{2p+1}_n$ (pour des permutations partitionnant $[n]$ en
un nombre impair d'orbites) plus un {\'e}l{\'e}ment de $\alpha M_c[n]+\delta
M_c[n]$. Ainsi, $\Phi_{\gl_0}(F_n)$ est inclus dans le sous-$\QQ$--espace
vectoriel de $\Sigma_n$ engendr{\'e} par les diagrammes ayant un nombre impair de
composantes connexes.
\item Enfin, on peut pr{\'e}ciser l'image de $\Phi_{\gl_0}(F_n)$ en utilisant le
  fait que, si sdim$(E)=0$, pour les morphismes canoniques
  $q\co \gl(E)\twoheadrightarrow \pgl(E)$ et $j\co \psl(E)\inj\pgl(E)$, on a:
  $q_*\Phi_{\gl(E)}(F_n)=j_*\Phi_{\psl(E)}(F_n)\subset q_*(\sll(E)^{\otimes
  n})$.\\ Soit $I$ l'{\'e}l{\'e}ment identit{\'e} de $\gl(E)$. Notons que $\Ker(q_*)$ est
  l'id{\'e}al engendr{\'e} par $I$ dans l'alg{\`e}bre tensorielle de $\gl(E)$. De plus
  $\sll(E)$ est l'orthogonal de $I$ dans $\gl(E)$ et ainsi l'orthogonal du noyau
  de $q_*$ est l'alg{\`e}bre tensorielle de $\sll(E)$. L'image de
  $\Phi_{\gl(E)}(F_n)$ est dans la somme de $\Ker(q_*)$ et de $\sll(E)^{\otimes
  n}$; l'orthogonal de cette somme est l'intersection de $\sll(E)^{\otimes n}$
  et de $\Ker(q_*)$, c'est l'espace engendr{\'e} par $\{x\otimes I\otimes
  y\in\sll(E)^{\otimes n}\}$.  Consid{\'e}rant $F_n\subset \D([n],[0])$ et en
  utilisant la dualit{\'e} dans $\gl(E)$, on a ainsi,\\
  $$\forall K\in F_n,\,\,\forall x\otimes I\otimes y\in\sll(E)^{\otimes
    n},\,\,\Phi_{\gl(E)}(K) (x\otimes I\otimes y)=0$$ 
  Si $i\leq n$, on note $<\sigma\setminus i>$ l'{\'e}l{\'e}ment $\Delta<\mu>$ si
  $\sigma(i)=i$ sinon l'{\'e}l{\'e}ment $<\mu>$ de $\Sigma_{n-1}$ o{\`u} $\mu$ est l'{\'e}l{\'e}ment
  de $\S_{n-1}$ qui, conjugu{\'e} avec la bijection croissante de $[n-1]$ vers
  $[n]\setminus\{i\}$, vaut $\mu'$ d{\'e}fini par $\mu'(j)=\sigma(j)$ si
  $\sigma(j)\neq i$, $\mu'(j)=\sigma(i)$ sinon.  Revenant {\`a} la d{\'e}finition de
  $\Phi_{gl(E),E}$, il est clair que si $x\in \gl(E)^{\otimes i-1}$ et $y\in
  \gl(E)^{\otimes n-i}$ alors
  $$\Phi_{gl(E),E}(<\sigma>) (x\otimes I\otimes
  y)=\Phi_{gl(E),E}(<\sigma\setminus i>) (x\otimes y)$$ 
  Enfin $\Phi_{gl(E),E}$ induit un morphisme de $\Sigma_n$ dans
  $(\gl(E)^{\otimes n})^*$ donc par restriction un morphisme de $\Sigma_n$ dans
  $(\sll(E)^{\otimes n})^*$. L'intersection des noyaux de ces morphismes
  (lorsque $E$ varie) est $\Sigma'_n$. En effet, si dim$(E)\geq n$, l'unique
  permutation de $\Sigma_1$ est envoy{\'e}e sur $(x\im \str_E(x))$ qui est nul sur
  $\sll(E)$ et donc toute permutation ayant un point fixe est envoy{\'e}e sur
  z{\'e}ro.\\
  R{\'e}ciproquement, on peut associer {\`a} la permutation $\sigma$ de $\S_n$ un tenseur
  $x_\sigma\in\gl(E)^{\otimes n}$ tel que $\Phi_{gl(E),E}(<\sigma>) (x_\mu)=1$
  si $\mu=\sigma$, 0 sinon. Pour construire de tels tenseurs, on identifie
  $\gl(E)$ {\`a} $\gl{(\dim(E_{\bar0}),\dim(E_{\bar1}))}(\QQ)$ et notant $e_{i,j}$
  les matrices {\'e}l{\'e}mentaires, on d{\'e}finit d'abord\\
  $x_{(1,2,\ldots k_1) (k_1+1,k_1+2,\ldots k_2)\ldots(k_p+1,\ldots
    n)}=(e_{1,2}\otimes e_{2,3}\otimes \ldots \otimes e_{k_1,1}) \otimes
  (e_{k_1+1,k_1+2}\otimes \ldots \otimes e_{k_2,k_1+1})\otimes \ldots \otimes
  (e_{k_p+1,k_p+2}\otimes \ldots \otimes e_{n,k_p+1})$\\ puis
  $x_{\mu\sigma\mu^{-1}}=\mu(x_\sigma)$; si $\sigma$ n'a pas de point fixe,
    $x_\sigma$ appartient {\`a} $\sll(E)^{\otimes n}$. Ainsi notant $\Sigma'_n$
  le sous-module de $\Sigma_n$ engendr{\'e} par les permutations ayant un point
  fixe, le morphisme suivant est nul:
  $$f_i\co F_n \stackrel{\Phi_{\gl_0}}{\go}{\Sigma_n}_{/(\Delta\Sigma_n)}
  \stackrel{<\sigma>\im <\sigma\setminus i>}{\go}
  {\Sigma_{n-1}}_{/(\Delta\Sigma_{n-1})}\go{\Sigma_{n-1}}_{/\Sigma'_{n-1}+
  (\Delta\Sigma_{n-1})}$$
\item On peut faire une derni{\`e}re remarque sur $\Phi_{\gl_0}(F_6)$: elle provient
  du fait que $\Phi_{\gl_0}$ est nul sur $F_3$ en degr{\'e} sup{\'e}rieur ou {\'e}gal {\`a} $3$.
  $\Phi_{\gl_0}(F_6)$ est donc inclus dans le noyau de l'application composition
  avec le diagramme suivant:
  $$\input{D4.tex}$$
\end{itemize}

\subsubsection{D{\'e}monstration de l'existence du caract{\`e}re $\X_3$}\label{demo X3}
Notons $\mu$ la forme lin{\'e}aire d{\'e}finie sur $F_6$ en degr{\'e} impair par
$\X_\sll(D_1\circ\tilde K)=-\mu\delta\beta^p$ modulo $(\delta^2)$ (avec les
notations d{\'e}j{\`a} utilis{\'e}es dans (3.4.2)). Nous allons montrer que $\mu=\eta$.\\
D'abord on a: $\Phi_\gl(D_1\circ\tilde K)\equiv
\mu(K)\Delta\left(<(1,2,3)>-<(1,3,2)>\right)$ modulo $\Delta^2\Sigma_{3}$. On
d{\'e}finit une forme lin{\'e}aire $\ g$ sur $\Sigma_{3}$ par
$g(\Sigma'_{3}+\Delta^2\Sigma_{3})=0$, $g(<\sigma>)=0$ si $\sigma$ n'est pas un
$3$--cycle, $g(<(1,2,3)>)=-g(<(1,3,2)>)=\frac12$. Soit $\mu_0\co \Sigma_6\go\QQ$
l'application qui associe {\`a} chaque {\'e}l{\'e}ment de $\Sigma_6$ l'image par $g$ de sa
compos{\'e}e avec $D_1$.  Alors $\mu(K)=\mu_0(\Phi_\gl(K))$ mais un calcul explicite
de l'application composition avec $D_1$ montre que l'image de cette application
est en fait dans $\Sigma'_{3}+\Delta\Sigma_{3}$, ainsi $\mu_0$ est nul sur
$\Delta\Sigma_6$. De plus, dans $\D_\gl$, l'{\'e}l{\'e}ment suivant est nul:
$${\epsfbox{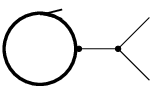}}$$
Donc $\mu_0(\Sigma'_6)=0$.\\ De m{\^e}me $\eta''$ induit une forme lin{\'e}aire
$\eta_0$ sur $\Sigma_6$ telle que
$\eta_0(<\sigma>)=\eta''(\Phi{\gl_{2,2}}(<\sigma>))$ pour $\sigma\in\Sigma_6$.
En particulier $\eta''$ est nulle sur l'orthogonal de $((\gl_{2,2})_{\bar
  1})^{\otimes 6}$ et cet orthogonal contient l'id{\'e}al Ker$(q)$ {\'e}voqu{\'e} dans la
question pr{\'e}c{\'e}dente. il en r{\'e}sulte que $\eta_0$ est nulle sur $\Sigma'_6$.\\ 

Des calculs effectu{\'e}s avec Maple donnent les r{\'e}sultats suivants:
\begin{itemize}
\item Le sous-espace de ${\Sigma_6}_{/\Sigma'_6}$ des {\'e}l{\'e}ments
$(G,\epsilon)$--invariants, invariants par l'endomorphisme de $\Sigma_6$ d{\'e}fini
sur sa base par $<\sigma>\im<\sigma>+(-1)^n<\sigma^{-1}>$, form{\'e} de combinaisons
lin{\'e}aires diagrammes ayant un nombre impair de composantes connexes, est de
dimension quatre, engendr{\'e} par les {\'e}l{\'e}ments suivants:\\
$y_1=f\left((1,2,3,4,5,6)\right)$\\ $y_2=f\left((1,2,3,5,4,6)\right)$\\
$y_3=f\left((1,3,2,5,4,6)\right)$\\ $y_4=f\left((1,3) (2,5) (4,6)\right)$\\
Avec $f(y_i)=\sum_{g\in G} \epsilon(g) (<gy_ig^{-1}>+<gy_i^{-1}g^{-1}>)$\\
L'image par $\Phi_{\gl_0}$ des {\'e}l{\'e}ments de degr{\'e} impair de $F_6$ est donc de
dimension inf{\'e}rieure ou {\'e}gale {\`a} quatre.
\item L'application $f_6$ (construite dans la section pr{\'e}c{\'e}dente {\`a} l'aide de
  l'applic\-ation $<\sigma>\im<\sigma\setminus 6>$) envoie $y_1$ et $y_4$ sur $0$,
  mais $f_6(y_2)=f_6(y_3)\neq 0$.
\item L'application induite par recollement de l'{\'e}l{\'e}ment $D_4$ dans $\D_\gl$
  envoie $y_1$, $y_2$, $y_3$ et $y_4$ sur respectivement $-24\Delta^2$,
  $16\Delta^2-16$, $-8\Delta^2-16$ et $-48$ fois l'{\'e}l{\'e}ment
  $(<(1,2,3)>-<(1,3,2)>)$ de $\Sigma_{3}$. En cons{\'e}quence, et {\`a} l'aide des deux
  derni{\`e}res remarques de la section pr{\'e}c{\'e}dente, on peut conclure que l'image par
  $\Phi_{\gl_0}$ des {\'e}l{\'e}ments de degr{\'e} impair de $F_6$ est de dimension inf{\'e}rieure
  o{\`u} {\'e}gale {\`a} deux, engendr{\'e}e par $y_1$ et $y_2-y_3$.
\end{itemize}
De plus, $\Phi_\gl(D_2)=\frac16y_1$ et $D_1\circ D_2=8t^3\in\L$ donc
$\mu_0(y_1)=6\mu(D_2)=0$. On va montrer que $\eta_0(y_1)$ est nul ce qui
permettra de conclure que $\eta_0=\mu_0$ sur l'image de $F_6$ en degr{\'e} impair
donc que $\eta=\mu$. En effet, sur cette image, $\eta_0$ et $\mu_0$ ont m{\^e}me
noyau donc sont proportionnelles, donc $\eta$ et $\mu$ sont proportionnelles, or
$\eta$ et $\mu$ co{\"\i}ncident et sont non nulles sur le diagramme suivant:\\
$${\epsfbox{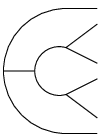}}$$
dont la compos{\'e}e avec $D_1$ donne l'{\'e}l{\'e}ment $x_3\in\L$; elles sont donc bien
{\'e}gales.\\
Pour montrer que $\eta_0(y_1)=0$, on calcule directement $\eta(D_2)$:\\
Le tenseur $x_1$ intervenant dans le calcul de $\eta$ est obtenu comme la
r{\'e}duction par l'application $\phi$ du diagramme colori{\'e} suivant:
$$\input{x1_calc.tex}$$
Donc $\eta(D_2)$ est obtenu comme la r{\'e}duction modulo  $(a)R+(b-1)R$ de:
$$\input{etaD2.tex}$$
La sommation {\'e}tant faite pour $\gamma_1$, $\gamma_2$ et $\gamma_3$ parcourant
$\{-a_1\omega_1,-a_2\omega_2,-a_3\omega_3\}$ (Si l'un des $\gamma_i$ valait
$\pi$, le diagramme colori{\'e} obtenu ne serait pas admissible). Or si $i\neq j$
alors
$$\input{phinul.tex}$$
donc tous les termes de la somme sont nuls sauf celui pour lequel
$\gamma_1=\gamma_2=\gamma_3=-a_1\omega_1$ mais ce terme est dans $a_1^2R$; ainsi
$\eta(D_2)=0$ donc $\eta=\mu$.
\medskip\\
Pour terminer, nous avons montr{\'e} l'{\'e}galit{\'e} du premier coefficient des caract{\`e}res
$\X_\DD(x)$ et $\X_\sll(x)$ lorsque $x$ est dans $(D_2)_*(F_6)$. Pour avoir
l'{\'e}galit{\'e} de ces coefficients sur $\Lambda$ tout entier, il suffit de remarquer
que si un diagramme de $\L$ n'est pas dans $(D_2)_*(F_6)$, alors il repr{\'e}sente
un {\'e}l{\'e}ment de $\L$ divisible par $t$ donc pour lequel les deux coefficients sont
nuls.  Cette derni{\`e}re remarque termine la d{\'e}monstration du lemme.

\subsection{Les caract{\`e}res exceptionnels et la relation du carr{\'e}}
\subsubsection{Th{\'e}or{\`e}me d'existence}\label{theo d'exis}
Si $L$ est l'une des cinq alg{\`e}bres de Lie exceptionnelles de la liste\\
$({\gg_2},{\ff_4},\ee_6,\ee_7,\ee_8)$ alors $I_L=(P_\ex)+(u-\alpha_L t^2)$, la
liste des valeurs de $\alpha_L$ est
$(\frac{7}{72},\frac{-4}{81},\frac{-22}{225},\frac{-7}{81},\frac{-5}{72})$; ces
valeurs sont distinctes deux {\`a} deux. On identifie le but de chacun de ces
caract{\`e}res $\X_L$ avec $S_{/I_L}\simeq\QQ[t]$. Posons $Q_\ex=\prod_L(u-\alpha_L
t^2)$ polyn{\^o}me de degr{\'e} dix de $S$.
\begin{lem}
  Il existe un caract{\`e}re $\X_\ex\co \L\go S_{/(P_\ex)+(Q_\ex)}$ factorisant les cinq
  caract{\`e}res exceptionnels.
\end{lem}
\begin{dem}
Pour factoriser $\X_{\gg_2}$ et $\X_{\ff_4}$ on remarque que $I_{\gg_2}\cap
I_{\ff_4}$ est l'id{\'e}al somme $(P_\ex)+((u-\alpha_1t^2)(u-\alpha_2t^2))$ (car
$S_{/P_\ex}\simeq\QQ[t,u]$ est factoriel) donc le diagramme commutatif suivant
est un carr{\'e} cart{\'e}sien:
\begin{eqnarray*}
  S_{/(P_\ex)+((u-\alpha_1t^2) (u-\alpha_2t^2))}&\go&S_{/(P_\ex)+(u-\alpha_1t^2)}\\ 
  \downarrow&&\downarrow\\ 
  S_{/(P_\ex)+(u-\alpha_2t^2)}&\go&S_{/(P_\ex)+(u-\alpha_1t^2)+(u-\alpha_2t^2)}
  \simeq \QQ[t]_{/(t^2)}
\end{eqnarray*}
Ainsi, il suffit de montrer que le diagramme suivant est commutatif:
\begin{eqnarray*}
  \L&\stackrel{\X_{\gg_2}}{\go}&S_{/(P_\ex)+(u-\alpha_1t^2)}\\ 
  \X_{\ff_4}\downarrow&&\downarrow\\ 
  S_{/(P_\ex)+(u-\alpha_2t^2)}&\go&S_{/(P_\ex)+(u-\alpha_1t^2)+(u-\alpha_1t^2)}
  \simeq \QQ[t]_{/(t^2)}
\end{eqnarray*}
Mais l'anneau en bas {\`a} droite est nilpotent et nul en degr{\'e} sup{\'e}rieur {\`a} deux. Or
$\L$ ne diff{\`e}re pas de $\L_0$ en degr{\'e} inf{\'e}rieur ou {\'e}gal {\`a} dix (cf \cite{Kn1}) et les
deux caract{\`e}res co{\"\i}ncident sur $\L_0$ comme l'indique la formule (\ref{somme des
xn}) donn{\'e}e section (\ref{Xslosp}). Ainsi il existe un premier caract{\`e}re
interm{\'e}diaire factorisant $\X_{\gg_2}$ et $\X_{\ff_4}$. On r{\'e}it{\`e}re le m{\^e}me
proc{\'e}d{\'e} (d'abord en rempla\c{c}ant $\X_{\gg_2}$ par le caract{\`e}re que l'on vient de
construire et $\X_{\ff_4}$ par $\X_{\ee_6}$ etc...) et le m{\^e}me argument permet
de construire {\'e}tape par {\'e}tape le caract{\`e}re $\X_\ex$.
\end{dem}
On d{\'e}signe par $K_1$, $K_2$ et $K_3$ les combinaisons lin{\'e}aires de diagrammes
suivantes:
\setcounter{equation}{0}
\begin{eqnarray}
&\begin{picture}(0,0)%
\epsfig{file=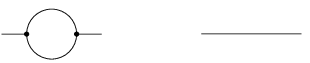}%
\end{picture}%
\setlength{\unitlength}{1579sp}%
\begingroup\makeatletter\ifx\SetFigFont\undefined%
\gdef\SetFigFont#1#2#3#4#5{%
  \reset@font\fontsize{#1}{#2pt}%
  \fontfamily{#3}\fontseries{#4}\fontshape{#5}%
  \selectfont}%
\fi\endgroup%
\begin{picture}(3621,625)(592,-371)
\put(2026,-136){\makebox(0,0)[lb]{\smash{\SetFigFont{12}{14.4}{\familydefault}{\mddefault}{\updefault}$-2t$}}}
\end{picture}
&\\
&\input{Rels3.tex}&\\
&\input{Rels2.tex}&
\end{eqnarray}
Soit $\L_1$ l'id{\'e}al de $\L\otimes S$ engendr{\'e} par les combinaisons lin{\'e}aires
de $([3],\emptyset)$--diagrammes connexes se factorisant par $K_1$ ou $K_2$ et
soit $\L_2$ l'id{\'e}al de $\L\otimes S$ engendr{\'e} par les combinaisons lin{\'e}aires
de $([3],\emptyset)$--diagrammes connexes se factorisant par $K_3\otimes K_3$.
On montrera le
\begin{lem}\label{ya des carres} L'application $S\go(\L\otimes
  S)_{/(\L_1+\L_2)}$ donn{\'e}e par l'unit{\'e} de $\L$ est surjective en degr{\'e}
  inf{\'e}rieur ou {\'e}gal {\`a} $20$. En cons{\'e}quence, il existe un caract{\`e}re $\X_4$ sur
  $\L$ en degr{\'e} inf{\'e}rieur ou {\'e}gal {\`a} $20$ {\`a} valeurs dans un quotient de $S$ qui
  factorise tous les caract{\`e}res annulant $K_1$, $K_2$ et $K_3\otimes K_3$.
\end{lem}
\paragraph{\bf Remarque}Des calculs {\'e}l{\'e}mentaires dans l'alg{\`e}bre
  d'endomorphisme\nl $\D([2],[2])\otimes S$ permettent de prouver que
  l'application $S\go(\L_0\otimes S)_{/(\L_1)}$ est surjective. Un ant{\'e}c{\'e}dent de
  $x_n$ peut {\^e}tre calcul{\'e} par la formule (\ref{somme des xn}) donn{\'e}e section
  \ref{Xslosp}.

On va aussi montrer que les caract{\`e}res $\X_3$ modulo $(P_\ex)+(R)$ et $\X_\ex$ mod
$(tP_{\sll}P_{\osp}P_{\sll_2})$ {\`a} valeurs dans l'anneau
$S_{/(tP_{\sll}P_{\osp}P_{\sll_2})+(P_\ex)+(R)}$ concentr{\'e} en degr{\'e} inf{\'e}rieur ou
{\'e}gal {\`a} $20$, se factorisent par $\X_4$, et donc co{\"\i}ncident, ce qui d{\'e}montrera
l'existence d'un caract{\`e}re $\X$ factorisant $\X_\ex$ et $\X_3$.\\
En fait, $\X_\ex$ annule $K_3$. Ceci provient du fait que pour chaque alg{\`e}bre
exceptionnelle, le carr{\'e} du Casimir engendre le sous-espace des {\'e}l{\'e}ments
$L$--invariants de $S^4L$ (Cf \cite{Vo2}).\\

\subsubsection{Unicit{\'e} de $\X_4$}
Pour d{\'e}montrer l'existence et l'unicit{\'e} de $\X_4$, il est plus ais{\'e} de
manipuler les diagrammes de $F_0$ qui forment un $\L$--module libre isomorphe {\`a}
$\L$. Il s'agit de montrer qu'en degr{\'e} inf{\'e}rieur ou {\'e}gal {\`a} $21$, $F_0$ est
isomorphe {\`a} l'espace $R_0$ engendr{\'e} par les diagrammes contenant l'un des trois
diagrammes suivants:
\begin{equation} \label{diag reductib}
\begin{array}{ccc}
\put(-20,-10){\epsfbox{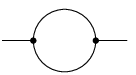}}&\put(-10,-10){\epsfbox{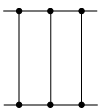}}
&\put(-40,-10){\epsfbox{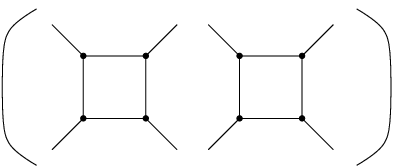}}\\
\qquad W_2\qquad&\qquad W_{2,2}\qquad&\qquad\qquad W_4\amalg W_4\qquad
\end{array}
\end{equation}
On introduit les notations:\label{NN}
$$N_3=\{(a,b,c)\in\NN^3 / a\leq b\leq c\}$$
muni de l'ordre lexicographique de $\NN^3$,
$$N_6=\{(\alpha,\beta)\in N_3\times N_3 / \alpha\leq \beta\}$$
muni de l'ordre induit par l'ordre lexicographique de $\NN^6$.
Dans la suite, on fixe $n\in\NN$, $\gamma=(\alpha,\beta)\in N_6$ et
$\delta=(a,b,c)\in N_3$. On d{\'e}finit l'ensemble
$$N=\NN\amalg N_3 \amalg N_6$$ que l'on finit d'ordonner en posant:\\
$\delta<\gamma$ si $\delta<\alpha$ et $\gamma<\delta$ si $\alpha\leq\delta$; si
$n<a$ alors $n<\delta$, sinon $\delta<n$;\\ $n$ et $\gamma$ sont ordonn{\'e}s comme
$n$ et $\alpha$;\\
l'ensemble $N$ est ainsi muni d'un bon ordre.\\ 
On notera $|\delta|=a+b+c$ et $|\gamma|=|\alpha|+|\beta|$.  Enfin, $W_\delta$ et
$W_n$ repr{\'e}senteront les diagrammes suivants de $F_{|\delta|}$ et $F_n$:
$${\epsfbox{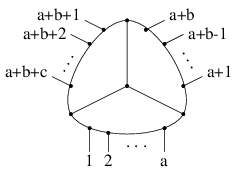}}W_\delta\qquad\qquad\qquad{\epsfbox{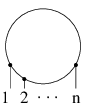}}\quad
W_n$$ 
et $W_\gamma$ le diagramme $W_\alpha\otimes W_\beta\in\D([0],[|\gamma|])$.  On
appelle ``roue'' un diagramme isomorphe {\`a} $W_n$ {\`a} la num{\'e}rotation des sommets
pr{\`e}s.\\ 
On d{\'e}signe par $f_K\co \D([n],[m])\go\D([0],[m])$ la composition {\`a} gauche par $K\in
\D([0],[n])$ (l'image de $f_K$ repr{\'e}sente les diagrammes ``contenant'' $K$) et
on pose:\\
$R^0_\delta=f_{W_\delta}($sous-espace de $\D([|\delta|],[0])$ engendr{\'e} par
les diagrammes connexes$)$\\
$R^0_\gamma=f_{W_\gamma}($sous-espace de $\D([|\alpha|]\amalg[|\beta|],[0])$
engendr{\'e} par les diagrammes dont chaque composante connexe a au moins un sommet
trivalent et rencontre $[|\alpha|]$ et $[|\beta|])$.\\
Ceci permet de d{\'e}finir pour $d\in\NN\setminus\{0\}$ et $\alpha\in N\setminus\{0\}$
$$R^0_d=f_{W_d}(F_d) \qquad R_\alpha^<=R_0+\sum_{\beta<\alpha}R^0_\beta$$
$$R_\alpha=R_\alpha^<+R^0_\alpha \qquad \hbox{et } \qquad \bar
R_\alpha=R_\alpha/R_\alpha^<$$
Si $\alpha<\beta$ dans $N$, on a $R_0\subset R_\alpha^<\subset R_\alpha\subset
R_\beta^<\subset R_\beta$.
\begin{lem}\label{decoroue}
  $F_n$ est engendr{\'e} en degr{\'e} $n$ par les roues et en degr{\'e} sup{\'e}rieur
  ou {\'e}gal {\`a} $n+1$ par les diagrammes $K=W\circ K'$ o{\`u} $K'$ est un arbre
  ({\'e}l{\'e}ment de $F_{k}$ de degr{\'e} $k-1$) vu comme {\'e}l{\'e}ment de $\D([0],[k])$
  et $W$ est une roue ({\`a} $n+k$ jambes). 
\end{lem}
\begin{dem}
La d{\'e}monstration de ce lemme est laiss{\'e}e au lecteur (elle
repose sur une simple manipulation des relations (IHX)).\\
{\bf Indications}\qua Il est utile de remarquer que l'on peut raisonner sur le
``squelette'' des diagrammes (c'est {\`a} dire oublier les jambes des
diagrammes qui peuvent {\^e}tre d{\'e}plac{\'e}es par des relations
(IHX)). L'{\'e}tape interm{\'e}diaire est de d{\'e}montrer que $F_n$
est engendr{\'e} en degr{\'e} sup{\'e}rieur ou {\'e}gal {\`a} $n+1$
par des diagrammes ayant un squelette de la forme suivante:
\vspace{.1cm}\\
$\put(20,0){\epsfbox{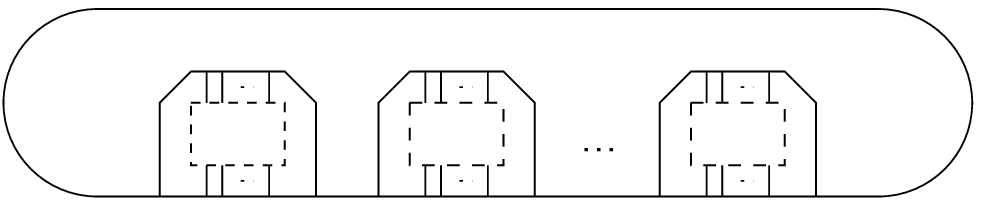}}\put(85,17){\small$\sigma_1$}\put(147,17)
{\small$\sigma_2$}\put(227,17){\small$\sigma_k$}$
\\ o{\`u} chacune des boites est donn{\'e}e par une permutation
$\sigma_i\in\S_{n_i}$.\\ 
Ensuite, il est possible de conclure grace {\`a} la manipulation de
``fusion des arbres'' suivante:
\vspace{.2cm}\\
$\put(10,-5){\epsfbox{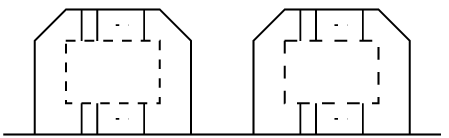}}\put(38,12){\small$\sigma_1$}\put(103,12)
{\small$\sigma_2$}\hspace{5cm}=\put(10,-5){\epsfbox{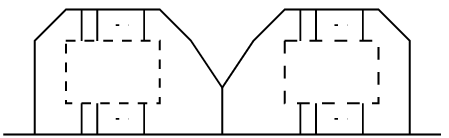}}\put(38,12){$\small\sigma_1$}\put(103,12)
{\small$\sigma_2$}\hspace{5cm}$\\
\hspace*{5cm}$+\put(10,-5){\epsfbox{fusionI.eps}}\put(38,12){\small$\sigma_1$}\put(111,15)
{\begin{rotate}{180}{\small$\sigma_2$}\end{rotate}}$\end{dem}

On utilise aussi le fait que $W_{2n+1}\in\Im(f_{W_{2n}})$ (en particulier pour
$n=2$) et de la m{\^e}me mani{\`e}re, l'image de $f_{W_{2,2}}$ contient le diagramme
suivant:
$$\epsfbox{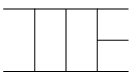}$$ 
Si $d\in\NN$, on notera $\D(d,[n])$ le quotient de $\D([d],[n])$ tuant l'action
du groupe $\S_d$. Si $\alpha\in N_3$, on note $\S_\alpha$ le sous-groupe de
$\S_{|\alpha|}$ isomorphe {\`a} $\S_{a_1}\times \S_{a_2}\times \S_{a_3}$ induit par
l'isomorphisme $[a_1]\amalg [a_2]\amalg [a_3]\inj [a_1+a_2+a_3]$ et
$\D(\alpha,[n])$ le module quotient de $\D([|\alpha|],[n])$ par les relations:
si $\sigma\in\S_\alpha$ et $K\in \D([|\alpha|],[n])$ alors $\sigma.K\equiv
K$. On a une application surjective naturelle de $\D([|\alpha|],n)$ vers
$\D(\alpha,n)$: celle ci revient {\`a} remplacer la num{\'e}rotation des sommets
univalents de la source par un coloriage {\`a} l'aide de trois couleurs que l'on
notera $x_1,\,x_2$ et $x_3$.\\
De m{\^e}me si $\gamma=(\alpha,\beta)\in N_6$, on note $\D(\gamma,[n])$ le quotient
de $\D([|\gamma|],[n])$ tuant l'action du groupe
$\S_\gamma\simeq\S_\alpha\times\S_\beta \subset\S_{|\gamma|}$.  L'application
quotient revient {\`a} remplacer la num{\'e}rotation des sommets univalents de la source
d'un diagramme par un coloriage {\`a} l'aide de six couleurs que l'on notera
$x_1,x_2,x_3$ et $y_1,y_2,y_3$.
\begin{lem}\label{sym} Si $\gamma\in N$ l'application 
  $$\bar f_{\gamma}\co \D([|\gamma|],[n])\go\bar R_\gamma$$ (obtenue par composition
de $f_{W_\gamma}$ avec l'application quotient $R_\gamma\go\bar R_\gamma$) se factorise
  par une application $g_\gamma\co \D(\gamma,[n])\go\bar R_\gamma$.
\end{lem}
Pour d{\'e}montrer ce lemme, il suffit de voir que pour tout
$K\in\D([|\gamma|],[n])$, pour toute transposition $\sigma=(i,i+1)\in
\S_\gamma$, on a $\bar f_\gamma(\sigma.K)=\bar f_\gamma(K)$. Mais
$W_\gamma\circ(\sigma.K)-W_\gamma\circ K=(\sigma.W_\gamma-W_\gamma)\circ K$ et
par la relation (IHX), $(\sigma.W_\gamma-W_\gamma)\circ K\in R^<_\gamma$.

\begin{lem} Soit $\gamma\in N$ et $K\in\D(\gamma,[0])$ alors $g_\gamma(K)$ est
nul ou $\dim(H_1(K))>22$. Ainsi, $F_0\subset R_0$ en degr{\'e} $d$ tel que $1<d<22$.
\end{lem}
\dem
Pour montrer ce dernier lemme, on proc{\`e}de par {\'e}tapes:
\begin{itemize}

\item Si $\gamma<((1,3,3),(3,3,3))$ il n'y a rien {\`a} d{\'e}montrer car
  $\Im(g_\gamma)\subset R_0$.
\item En utilisant les arguments du lemme \ref{sym}, il est facile de montrer
  qu'un diagramme contenant une roue {\`a} $d$ jambes appartient {\`a} $R_d$. En
  effet, si le compl{\'e}mentaire de la roue n'est pas connexe, quitte {\`a} permuter
  ses jambes, on peut se ramener {\`a} un diagramme obtenu en recollant deux
  {\'e}l{\'e}ments de $F_2$ et donc divisible par $t$ (i.e. appartenant {\`a} $R_0$).
\item Si $d\geq4$, par le lemme \ref{sym} et en appliquant le lemme
  \ref{decoroue} au compl{\'e}men\-taire de $W_d$ dans un diagramme de $R_d$, on
  obtient facilement $\bar R_d=0$.
\item Si les lettres $a$ et $b$ d{\'e}signent les deux lettres $\{x,y\}$, si
  $\{i,j,k\}\subset\{1,2,3\}$ et si $K$ contient l'un des cinq diagrammes
  suivants: 
  \begin{equation}\label{red1}\epsfbox{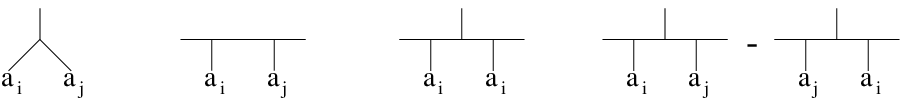}\end{equation} 
  alors $g_\gamma(K)=0$. En effet, on peut toujours faire appara{\^\i}tre
  une roue {\`a} moins de cinq jambes en recollant un tel diagramme {\`a}
  $W_\gamma$. 
\item On peut maintenant montrer que si $\alpha\in N_3$ alors $\bar R_\alpha$
est nul. Pour cela, appliquons le lemme (\ref{decoroue}) a un diagramme connexe
$K\in \D(\alpha,[0])$. Si $K$ est un arbre, il contient le premier diagramme de
(\ref{red1}) sinon, il se d{\'e}compose en une roue dont certaines jambes sont
colori{\'e}es par les $x_i$ et d'autres qui sont reli{\'e}es aux feuilles de
l'arbre. Aux moins trois d'entre elles sont de couleur $x_2$. S'il ne contient
pas le deuxi{\`e}me diagramme de (\ref{red1}), les jambes colori{\'e}es par les $x_i$
sont s{\'e}par{\'e}es par des feuilles de l'arbre. Nous dirons que deux jambes colori{\'e}es
sont en relation si elles ne sont s{\'e}par{\'e}es que par une feuille de l'arbre. On
compl{\`e}te cette relation en une relation d'{\'e}quivalence. Le dernier {\'e}l{\'e}ment de
(\ref{red1}) prouve que l'on peut permuter les couleurs des jambes qui sont en
relation sans modifier l'image par $g_\alpha$ de $K$. Mais alors le troisi{\`e}me
diagramme de (\ref{red1}) {\'e}tant annul{\'e} par $g_\alpha$, si $g_\alpha(K)\neq0$,
aucune des classes d'{\'e}quivalences ne peut contenir deux jambes colori{\'e}es par la
m{\^e}me couleur et il y a donc au moins trois classes d'{\'e}quivalences. Ceci signifie
qu'en au moins trois endroits dans $K$, deux feuilles de l'arbre sont reli{\'e}es
par une ar{\^e}te de la roue. La r{\'e}union de ces trois ar{\^e}tes et de l'arbre forme un
diagramme {\`a} trois boucles qui, par le lemme (\ref{decoroue}), permet de faire
appara{\^\i}tre un diagramme $W_\beta$ o{\`u} $\beta\in N_3$. Ainsi $f_{W_\alpha}(K)$
appartient {\`a} $R_{(\alpha,\beta)}\subset R^<_\alpha$.
\item $$\input{red_diah.tex}$$  
  de plus ces {\'e}l{\'e}ments sont nuls si $i$=$k$.
\item \begin{equation}\label{red2}\input{red_diai.tex}\end{equation}  
  de plus ces {\'e}l{\'e}ments sont nuls si $n\geq 3$.
\item Maintenant si $\gamma\in N_6$ et si $K\in\D(\gamma,[0])$ est un diagramme
  dont l'une des composantes connexes est un arbre, alors $K$ est combinaison
  lin{\'e}aire de diagrammes contenant l'un des deux premiers {\'e}l{\'e}ments de (\ref{red1}).
\item Si une composante connexe $K_0$ de $K$ est une roue, alors $K_0$ contient
  un diagramme de type (\ref{red2}) et soit $K_0\in R_4$, soit le nombre $n$ est
  sup{\'e}rieur ou {\'e}gal {\`a} quatre donc $g_\gamma(K)=0$.  
\item Si $K$ est connexe, on lui applique le lemme (\ref{decoroue}). Reprenant
la relation d'{\'e}quivalence pour les jambes de la roue colori{\'e}es par les $x_i$, on
retrouve qu'en trois endroits distincts, la roue poss{\`e}de deux jambes
cons{\'e}cutives qui ne sont pas colori{\'e}es par $x_i$ et au moins une des deux est
donc une feuille de l'arbre. L'arbre a donc au moins trois feuilles et $K$ a
donc au moins trois boucles.
\item Ainsi si $g_\gamma(K)\neq0$, soit $K$ est connexe et a au moins trois
boucles, soit il poss{\`e}de plusieurs composantes connexes ayant chacune au moins
deux boucles. Dans les deux cas, le nombre de sommets trivalents de $K$ est
sup{\'e}rieur ou {\'e}gal {\`a} $|\gamma|+4$. Comme $\gamma\geq ((1,3,3),(3,3,3))$, on a
$|\gamma|\geq 16$ et le nombre de sommets trivalents de $f_{W_\gamma}(K)$ est
sup{\'e}rieur ou {\'e}gal {\`a} $44$.\endproof
\end{itemize}

On peut maintenant montrer le lemme \ref{ya des carres}. En effet, consid{\'e}rons
un diagramme $K=K'\circ K''\in F_0$ o{\`u} $K'$ est l'un des trois {\'e}l{\'e}ments de
(\ref{diag reductib}). Le lecteur pourra v{\'e}rifier que si $K''$ n'est pas
connexe, alors, si $K'$ est isomorphe {\`a} $W_2$, $K$ est nul, si $K'$ est
isomorphe {\`a} $W_{2,2}$, $t$ divise $K$ et si $K'$ est la r{\'e}union disjointe de
deux carr{\'e}s, alors $K$ se d{\'e}compose en une combinaison lin{\'e}aire de diagrammes du
m{\^e}me type pour lesquels $K''$ est connexe ou $K'$ est l'un des deux premiers
{\'e}l{\'e}ments de (\ref{diag reductib}). On peut ainsi toujours se ramener au cas o{\`u}
$K''$ est connexe. Alors, $K$ se d{\'e}compose dans $(\L\otimes S)_{/(\L_1+\L_2)}$
en une combinaison lin{\'e}aire de diagrammes de degr{\'e} strictement inf{\'e}rieur et en
r{\'e}it{\'e}rant le processus pour un diagramme de $(\L\otimes S)_{/(\L_1+\L_2)}$ de
degr{\'e} inf{\'e}rieur ou {\'e}gal {\`a} $21$ (correspondant {\`a} un {\'e}l{\'e}ment de $\L$ de degr{\'e}
inf{\'e}rieur ou {\'e}gal {\`a} $20$), on peut l'exprimer comme combinaison lin{\'e}aire
d'{\'e}l{\'e}ments de $S$.
                                %
\subsubsection{$\X_3$ et les relations exceptionnelles}\label{X3ex}
Le but de cette section est de montrer que $\X_3$ modulo l'id{\'e}al engendr{\'e} par
$P_\ex$ annule $\L_2$. Les variables $t$, $u$ et $v$ sont d{\'e}termin{\'e}es par le
fait que chaque superalg{\`e}bre de Lie simple annule les {\'e}l{\'e}ments $K_1$ et $K_2$
d{\'e}finis dans la section \ref{theo d'exis}. Ainsi les caract{\`e}res $\X_3$ et
$\X_\ex$ annulent $\L_1$.\\
Le but de $\X_3$ modulo $(P_\ex)$ est l'anneau
\begin{eqnarray*}
&S_{/(P_\ex)+(tP_{\sll}P_{\osp}P_{\sll_2})}=
S_{/(P_\ex)+((9u-2t^2)(u-t^2)(9u-5t^2)t^3u^2)}\\
\inj&S_{/I_{\sll_3}}\times S_{/I_{\osp(1,2)}}\times S_{/(P_{\sll2})+(9u-5t^2)}
\times S_{/(P_{\ex})+(u^2)}\times S_{/(P_{\ex})+(t^3)}
\end{eqnarray*}
Tout comme les alg{\`e}bres de Lie exceptionnelles, pour $L=\sll_3$ et $L=\osp(1,2)$,
le carr{\'e} du Casimir engendre le sous-espace $L$--invariant de $S^4L$. Ainsi les
caract{\`e}res $\X_{\sll_3}$ et $\X_{\osp(1,2)}$ annulent $\L_2$. De plus,
$\Phi_{\sll_2}$ annule l'{\'e}l{\'e}ment
\begin{equation}\label{relsl2}
\epsfbox{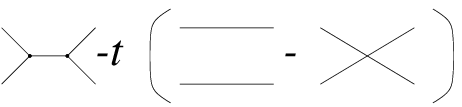}
\end{equation}
et ceci suffit {\`a} d{\'e}terminer le caract{\`e}re $\X_{\sll_2}$. Ainsi on peut v{\'e}rifier
directement que $(9u-5t^2)$ divise $\Phi_{\sll_2}(K_3)$.\\
De m{\^e}me $(\delta-8)$ divise $\Phi_\osp(K_3)$ et par homog{\'e}n{\'e}it{\'e} de $\X_\osp$, on
a: $\X_\osp(\L_2)\subset (\delta-8\alpha)^2\QQ[\alpha,\delta]$. Ceci montre
que $\X_3$ modulo $(P_{\ex})+(u^2)$ est nul sur $\L_2$.\\
Pour montrer que $\X_3$ modulo $(P_\ex)+(t^3)$ annule $\L_2$, il suffit de
montrer que les images de $\L_2$ par respectivement $\X_\sll$ et $\X_\DD$ sont
dans les id{\'e}aux respectifs $(t^3)$ et $(v^3)$. Pour cela, supposons qu'un
{\'e}l{\'e}ment $K$ de $F_0\otimes S$ soit de la forme suivante: 
$$\input{car2t3.tex}$$ o{\`u} $K'\circ K''$ est connexe. On note $u$ l'{\'e}l{\'e}ment
de $\L$ correspondant {\`a} $K$. Il se trouve que $\Delta$ divise $\Phi_\sll(K_3)$
et $a_1$ divise $\Phi_\DD(K_3)$ (ceci provient du fait que
$\Phi_{\psl(2,2)}(K_3)=0$).\\
Comme $K'\circ K''$ est connexe, $K'$ appartient {\`a} $\D_b([8],[6])$ et donc $a_1^2$
divise $\Phi_\DD((K_3\otimes K_3)\circ K')$ or il existe une forme lin{\'e}aire sur
$X^{\otimes 6}$ {\`a} valeurs dans $R$ (cf section \ref{prop D21}) qui prend sur
$a_1^{-2}\Phi_\DD((K_3\otimes K_3)\circ K')$ la valeur
$\frac1{a_1^3}\X_\DD(u)$. Ainsi $v^3$ divise $\X_\DD(u)$.\\
De m{\^e}me, en notant $\D'_\gl([p],[q])$ le sous-$\QQ[\Delta]$--espace de
$\D_\gl([p],[q])$ engendr{\'e} par les {\'e}l{\'e}ments de $\Sigma'_{p+q}$ on peut ais{\'e}ment
v{\'e}rifier que la composition {\`a} droite par un {\'e}l{\'e}ment de $\D_b$ laisse stable
$\D'_\gl$ (il suffit de le v{\'e}rifier pour des diagrammes de la forme
$[.,.]\otimes Id$). Ainsi, comme $\Phi_\gl(K_3\otimes K_3)\in
\Delta^2\Sigma_8+\Sigma'_8$, on a $\Phi_\gl((K_3\otimes K_3)\circ
K')\in\Delta^2\Sigma_6+\Sigma'_6$ et donc $\Delta^2$ divise
$\Phi_\sll((K_3\otimes K_3)\circ K')$. Or il existe une forme lin{\'e}aire sur
$\Sigma_6$, nulle sur $\Sigma'_6$, {\`a} valeur dans $\QQ$ (cf la forme $\mu_0$ de
la section \ref{demo X3}) qui prend sur $\Delta^{-2}\Phi_\sll((K_3\otimes
K_3)\circ K')$ la valeur $(\frac1{t^3}\X_\sll(u))$ modulo $(t)$. Ainsi, $t^3$
divise $\X_\sll(u)$.\\

Donc $\X_3(u)$ est bien dans l'id{\'e}al somme $(P_\ex)+(t^3)$. Le fait que les
diagrammes du type de $K$ engendrent $\L_2$ modulo $\L_1$ r{\'e}sulte de la remarque
faite {\`a} la fin de la section pr{\'e}c{\'e}dente. Ainsi, $\X_3$ modulo $(P_\ex)$ annule
aussi $\L_2$ et par suite se factorise par $\X_4$ en degr{\'e} inf{\'e}rieur ou {\'e}gal {\`a}
$20$. Ceci termine la d{\'e}monstration de l'existence de $\X$.

\section{Les cas de $\g(3)$, $\ff(4)$ et les branchements}
Le but de cette section est de d{\'e}montrer que le caract{\`e}re $\X$ factorise aussi
les caract{\`e}res $\X_{\g(3)}$ et $\X_{\ff(4)}$.\\
En fait, on montre que $\X_{\g(3)}$ se factorise par $\X_{\sll_2}$ et que
$\X_{\ff(4)}$ et $\X_{\sll_3}$ co{\"\i}ncident. Pour calculer $\X_{\g(3)}$, il suffit
de remarquer que le sous-module $X_2$ de $\L^2\g(3)$ form{\'e} par le noyau du
Casimir est simple et de superdimension nulle (pour $\sll_2$ ce m{\^e}me module est
nul). Notons $K_0$ l'{\'e}l{\'e}ment (\ref{relsl2}) annul{\'e} par $\Phi_{\sll_2}$ (cf
section \ref{X3ex}). Les {\'e}l{\'e}ments de $F_0$ de la forme
$$\begin{picture}(0,0)%
\epsfig{file=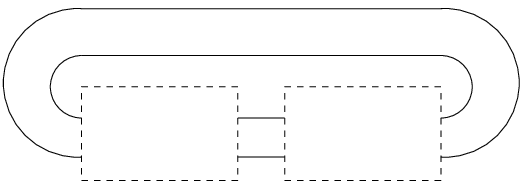}%
\end{picture}%
\setlength{\unitlength}{1973sp}%
\begingroup\makeatletter\ifx\SetFigFont\undefined%
\gdef\SetFigFont#1#2#3#4#5{%
  \reset@font\fontsize{#1}{#2pt}%
  \fontfamily{#3}\fontseries{#4}\fontshape{#5}%
  \selectfont}%
\fi\endgroup%
\begin{picture}(4969,1674)(591,-1273)
\put(4051,-886){\makebox(0,0)[b]{\smash{\SetFigFont{12}{14.4}{\rmdefault}{\mddefault}{\updefault}\small$K$}}}
\put(2101,-886){\makebox(0,0)[b]{\smash{\SetFigFont{12}{14.4}{\rmdefault}{\mddefault}{\updefault}\small$K_0$}}}
\end{picture}
$$ 
o{\`u} $K$ est un diagramme connexe de $\D([2],[2])$, sont
envoy{\'e}s par $\Phi_{\g(3)}$ sur $t\,\str_{X_2}(\Phi_{\g(3)}(K))=0$ et par suite
$\X_{\g(3)}=\X_{\sll_2}$ sont tous deux d{\'e}termin{\'e}s de mani{\`e}re unique par le fait
qu'ils annulent ces {\'e}l{\'e}ments.
\medskip\\

Supposons maintenant qu'une superalg{\`e}bre de Lie $L$ munie d'un {\'e}l{\'e}ment de
Casimir non d{\'e}g{\'e}n{\'e}r{\'e} $\Omega\in L\otimes L$ contienne une sous-alg{\`e}bre de Lie
$l$ sur laquelle la forme bilin{\'e}aire supersym{\'e}trique de $L$ ne soit pas
d{\'e}g{\'e}n{\'e}r{\'e}e. On a alors un foncteur $F\co Mod_L\go Mod_l$ qui consiste {\`a} regarder la
structure de $l$--module d'un $L$--module.\\
Si de plus le $l$--module $E$, orthogonal de $l$ dans $L$, v{\'e}rifie
$[E,E]_L\subset l$, on dira que $(l,E)$ est une bonne d{\'e}composition de $L$. On
peut remarquer que le Casimir de $L$ se d{\'e}compose en $\Omega=\omega+\pi$ avec
$\omega\in l\otimes l$ et $\pi\in E\otimes E$.\\
Pour transcrire cette situation en termes de diagrammes, on pose les d{\'e}finitions
suivantes:\\
Un $(X_1,X_2)$--diagramme bicolore est la donn{\'e}e d'un $(\Gamma,X_1)$--diagramme o{\`u}
$\Gamma$ est non orient{\'e}e et la donn{\'e}e d'un isomorphisme: $\partial\Gamma\simeq
X_2$. Nous dirons qu'une ar{\^e}te d'un diagramme bicolore est de la premi{\`e}re
couleur si elle n'appartient pas {\`a} la courbe $\Gamma$, et nous dirons qu'elle
est de la deuxi{\`e}me couleur dans le cas contraire. On repr{\'e}sentera toujours d'un
trait gras les ar{\^e}tes de la deuxi{\`e}me couleur. On note $\wh\A(X_1,X_2)$ le
$\QQ$--espace vectoriel de base les $(X_1,X_2)$--diagrammes bicolores quotient{\'e}
par les relations $(AS)$, $(IHX)$, $(STU)$. On note aussi
$$\wh\A(X)=\bigoplus_{X_1\amalg X_2=X} \wh\A(X_1,X_2).$$
On remarque que si $\Gamma$ est une courbe sans bord non orient{\'e}e,
$\A(\Gamma,X)\subset\wh\A(X,\emptyset)$.\\

Enfin on d{\'e}signe par $\wb\A(X_1,X_2))$ le quotient de $\wh\A(X_1,X_2)$ par les
relations not{\'e}es $(\wb{IHX})$:
$$\begin{picture}(0,0)%
\epsfig{file=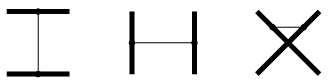}%
\end{picture}%
\setlength{\unitlength}{1973sp}%
\begingroup\makeatletter\ifx\SetFigFont\undefined%
\gdef\SetFigFont#1#2#3#4#5{%
  \reset@font\fontsize{#1}{#2pt}%
  \fontfamily{#3}\fontseries{#4}\fontshape{#5}%
  \selectfont}%
\fi\endgroup%
\begin{picture}(3088,688)(857,-405)
\put(3001,-161){\makebox(0,0)[b]{\smash{\SetFigFont{9}{10.8}{\rmdefault}{\mddefault}{\updefault}$-$}}}
\put(1801,-161){\makebox(0,0)[b]{\smash{\SetFigFont{9}{10.8}{\rmdefault}{\mddefault}{\updefault}$\equiv$}}}
\end{picture}
$$
et $\wb\A(X)=\bigoplus_{X_1\amalg X_2=X} \wb\A(X_1,X_2)$.

On d{\'e}finit les cat{\'e}gories $\wh{\D}$ et $\wb{\D}$ de mani{\`e}re analogue {\`a} $\D$
comme les cat{\'e}gories ayant les m{\^e}mes objets que $\D$ et dont les morphismes
sont\\ \label{Dbic}
$\wh{\D}([p],[q])=\wh\A([p]\amalg[q])$ et
$\wb{\D}([p],[q])=\wb\A([p]\amalg[q])$.  La composition dans $\wh{\D}$ ou
$\wb\D$ de deux diagrammes est encore obtenue par leur recollement si les ar{\^e}tes
issues des sommets de la source du premier diagramme sont de la m{\^e}me couleur
que celles issues du but du deuxi{\`e}me diagramme et on d{\'e}cr{\`e}te que la composition
est nulle dans le cas contraire. On note $\wb D\co \wh\D\go\wb\D$ le foncteur
quotient.

Soit $l$ une superalg{\`e}bre de Lie quadratique et $E$ un $l-module$ muni d'une
forme bilin{\'e}aire supersym{\'e}trique $l$--invariante. On note $\pi$ l'{\'e}l{\'e}ment de
$S^2E$ associ{\'e} et $\omega$ le Casimir de $l$.\\
\begin{prop}\label{moddim0}
  Il existe un unique foncteur mono{\"\i}dal $\QQ$--lin{\'e}aire
  $$\wh\Phi_{l,E}\co \wh{\D}\go\Mod_l$$ envoyant $[1]$ sur $l\oplus E$,
  prenant les m{\^e}mes valeurs que $\Phi_{l,E}$ sur $\A(\Gamma,X)$ lorsque $\Gamma$
  est une courbe sans bord non orient{\'e}e, et v{\'e}rifiant
  $$\input{foncth.tex}$$

  Si de plus $(l,E)$ est une bonne d{\'e}composition de $L$, alors $\wh\Phi_{l,E}$
  passe au quotient par $\wb D$, d{\'e}finissant un foncteur
  $\wb\Phi_{l,E}\co \wb\D\go\Mod_l$ qui v{\'e}rifie
  $$\wh\Phi_{l,E}=\wb\Phi_{l,E}\circ\wb D$$        
  Il existe un foncteur mono{\"\i}dal $\QQ$--lin{\'e}aire $\wb\Phi\co \D\go\wb\D$ d{\'e}fini 
  de mani{\`e}re unique par ses valeurs sur les morphismes suivants:
  $$\input{fonct.tex}$$ 
  Enfin, le foncteur $\wb\Phi$ v{\'e}rifie:
  $$F\circ\Phi_L=\wb\Phi_{l,E}\circ \wb\Phi$$
\end{prop}
\begin{dem}
  Nous justifierons seulement les existences de $\wb\Phi$ et de $\wbP$.  Notons
  provisoirement $f$ l'application qui associe {\`a} un $(X_1,X_2)$--diagramme
  bicolore le $(X_1\amalg X_2)$--diagramme sous-jacent (on oublie l'information
  sur les ``couleurs''). Ainsi il n'est pas difficile de voir que pour un
  diagramme $K\in\A(\emptyset,X)$ on a: $\wb\Phi(K)=\sum_{f(\wb K)=K}\wb K$. Le
  point cl{\'e} de la validit{\'e} de cette d{\'e}finition est qu'en notant $I$, $H$ et $X$
  les diagrammes intervenant dans la relation $(IHX)$, on peut r{\'e}organiser la
  somme $\sum_{f(\wb K)=I}\wb K - \sum_{f(\wb K)=H}\wb K + \sum_{f(\wb K)=X}\wb
  K$ de mani{\`e}re {\`a} faire appara{\^\i}tre la relation $(IHX)$ plus des relations
  $(STU)$ plus la relation $(\wb{IHX})$. En ce sens, les relations $(\wb{IHX})$
  sont n{\'e}cessaires et suffisantes {\`a} l'existence de $\wb\Phi$.\\
  En utilisant $\Omega=\omega+\pi$ et en ``d{\'e}veloppant'' le calcul de
  $\Phi_L(K)$ comme il a {\'e}t{\'e} fait dans la d{\'e}monstration du lemme \ref{PhiDD}
  pour l'alg{\`e}bre $\DD$, la formule $F\circ\Phi_L=\wbP\circ \wb\Phi$ appara{\^\i}t comme un
  simple jeu d'{\'e}criture et justifie du m{\^e}me coup l'existence de $\wbP$: en
  effet en ``d{\'e}veloppant'' $\Phi_L(I-H+X)=0$ on obtient bien que $\wh\Phi_{l,E}$
  v{\'e}rifie la relation $(\wb{IHX})$.
\end{dem}
{\bf Corollaire}\qua \sl Si $E=\bigoplus_i E_i$ o{\`u} chaque $E_i$ est un $l$--module de
superdimension nulle et tel que End$_l(E_i)\simeq\QQ$ alors les restrictions de
$\Phi_L$ et $\Phi_l$ {\`a} $F_0$ sont des formes lin{\'e}aires {\'e}gales.\rm\\ 
\begin{dem}
  En effet, l'image par $\wb\Phi$ d'un {\'e}l{\'e}ment de $F_0$ est {\'e}gale au m{\^e}me {\'e}l{\'e}ment vu
  dans $\wb\D([0],[0])$ plus une combinaison lin{\'e}aire de
  $(\Gamma,\emptyset)$--diagrammes o{\`u} $\Gamma\neq\emptyset$. Mais ces derniers
  s'interpr{\`e}tent comme la supertrace sur $E$ d'un tenseur $l$--invariant et sont
  annul{\'e}s par $\wb\Phi_{l,E}$ sous les hypoth{\`e}ses du corollaire. On a donc dans ces
  conditions $\wb\Phi_{l,E}\circ \wb\Phi=\Phi_l$ sur $F_0$.
\end{dem}
\begin{prop} La superalg{\`e}bre de lie $\ff(4)$ satisfait aux conditions ci-dessus
  pour $l=\sll(4,1)$ et $E$ est alors un $\sll(4,1)$--module simple de
  superdimension nulle. En cons{\'e}quence, les caract{\`e}res $\X_{\sll(4,1)}$ et
  $\X_{\ff(4)}$ sont {\'e}gaux.
\end{prop}
\begin{dem}\\
La partie paire de $\ff(4)$ est isomorphe {\`a} l'alg{\`e}bre semi-simple $\sll_2\times
\so_7$. Sa partie impaire est isomorphe au produit tensoriel de la
repr{\'e}sent-\break ation standard de $\sll_2$ par la repr{\'e}sentation spin$_7$. Consid{\'e}rons
une d{\'e}compos\-ition de Cartan: $\sll_2=\CC H\oplus \CC E\oplus \CC F$ avec
$[H,E]=2E$, $[H,F]=-2F$ et $[E,F]=H$. On note $V_2$ la repr{\'e}sentation standard
de $\sll_2$, $V_4$ la repr{\'e}sentation standard de $\sll_4$, $V'_4$ sa
repr{\'e}sentation duale, $W=\L^2V_4\simeq\L^2V'_4$. On choisit $e$ vecteur de plus
haut poids de $V_2$ et $f=F.e$.\\ $W$ est un $\sll_4$--module simple de dimension
$6$ autodual. Le choix d'une base de $W$ donne un mophisme d'alg{\`e}bre de Lie
$\sll_4\go\so_6$ qui est en fait un isomorphisme.\\ Fixons une injection
$\sll_4\simeq\so_6\inj\so_7$ et consid{\'e}rons la d{\'e}composition de $\ff(4)$ comme
$\CC H\times\sll_4$--module:
\begin{eqnarray*}
\ff(4)\simeq&\sll_2\oplus\so_7&\oplus V_2\otimes\spin_7\\
\simeq&\CC H\oplus \CC E\oplus\CC F\oplus \sll_4\oplus W&\oplus e\otimes V\oplus
e\otimes V'\oplus f\otimes V\oplus f\otimes V'
\end{eqnarray*}
Le crochet de Lie de $\ff(4)$ est un morphisme de $\CC H\times\sll_4$--module et
on remplit facilement la table du crochet suivante en utilisant la propri{\'e}t{\'e}
{\'e}tablie par V. G. Kac (cf \cite{Kac}) que pour toute superalg{\`e}bre de Lie classique
basique $\g$, si $\Delta$ est l'ensemble de ses racines, et si $\g_\alpha$
d{\'e}signe l'espace propre associ{\'e} {\`a} la racine $\alpha$, et $\beta\in\Delta$ est
diff{\'e}rente de $-\alpha$ alors:
\begin{equation}\label{propKac}
[\g_\alpha,\g_\beta]\neq 0\Longleftrightarrow \alpha+\beta\in\Delta
\end{equation}
$$\begin{array}{|c||c|c|c|c|c|c|c|c|c|} 
\hline
[.,.] & \CC H & \sll_4 & \CC E & \CC F & W & e\otimes V & f\otimes V' & e\otimes
V' & f\otimes V\\
\hline
\hline
\CC H & 0 & 0 & \CC E & \CC F & 0 & e\otimes V & f\otimes V' & e\otimes V' &
f\otimes V\\
\hline
\sll_4 & & \sll_4 & 0 & 0 & W & e\otimes V & f\otimes V' & e\otimes V' &
f\otimes V\\
\hline
\CC E  & & & 0 & \CC H & 0 & 0 & e\otimes V' & 0 & e\otimes V\\ 
\hline
\CC F & & & & 0 & 0 & f\otimes V & 0 & f\otimes V' & 0\\ 
\hline
W & & & & & \sll_4 & e\otimes V' & f\otimes V & e\otimes V & f\otimes V'\\
\hline
e\otimes V & & & & & & 0 & (\CC H)\oplus\sll_4 & \CC E & W\\ 
\hline
f\otimes V' & & & & & & & 0 & W & \CC F\\
\hline
e\otimes V' & & & & & & & & 0 & (\CC H)\oplus\sll_4\\ 
\hline
f\otimes V & & & & & & & & & 0\\
\hline
\end{array}$$
Le remplissage de cette table d{\'e}coule directement de la propri{\'e}t{\'e} de
surjectivit{\'e} du crochet signal{\'e}e ci-dessus {\`a} l'exception des termes $(\CC H)$
not{\'e}s entre parenth{\`e}ses pour lesquels on a par exemple:
$$[[e\otimes V,f\otimes V'],E]= [e\otimes V,[f\otimes V',E]]= [e\otimes
V,e\otimes V']= \CC E$$
Ce qui prouve que $\CC H\subset[e\otimes V,f\otimes V']$.\\

Posons maintenant
$$l=\CC H \oplus \sll_4 \oplus (e\otimes V) \oplus (f\otimes V')$$
$$X=\CC E \oplus \CC F \oplus W \oplus (e\otimes V') \oplus (f\otimes V)$$
On lit facilement sur la table que:
\begin{itemize}
\item $l$ est une sous-alg{\`e}bre de Lie de $\ff(4)$.
\item L'id{\'e}al engendr{\'e} par n'importe lequel de ses {\'e}l{\'e}ments non nul est $l$ tout
entier. Ainsi $l$ est simple et la classification de \cite{Kac} permet d'identifier
$l\simeq\sll(4,1)$.
\item $X$ est un $l$--module simple car il est monog{\`e}ne, engendr{\'e} par n'importe
lequel de ses {\'e}l{\'e}ments. 
\item Dans $\ff(4)$, $[X,X]$ est inclus dans $l$.
\end{itemize}
Compte tenu que pour une superalg{\`e}bre de Lie quadratique classique $\g$, si
$\alpha$ et $\beta$ sont des racines, $\g_\alpha\bot\g_\beta$ si et seulement si
$\alpha+\beta\neq0$, on a $l\bot X$. De plus, la superdimension de $X$ est bien
nulle comme annonc{\'e} dans la proposition.
\end{dem}
{\bf Remarque}\qua Le corollaire de la proposition \ref{moddim0} permet de
red{\'e}montrer que les caract{\`e}res $\X_{\sll(E)}$ et $\X_{\osp(E)}$ ne d{\'e}pendent que
de la superdimension de $E$, d'o{\`u} $\X_{\ff(4)}=\X_{\sll_3}$.

\renewcommand\refname{R\'ef\'erences}

\vfill

\end{document}